\documentclass{article}

\usepackage{appendix}
\usepackage{graphicx}       
\usepackage{amsmath}
\usepackage{amsfonts}
\usepackage{amssymb}
\usepackage{hyperref}

\usepackage{amsthm}

\usepackage{color}
\usepackage{graphicx}
\usepackage{mathtools, balance}
\usepackage{xcolor} %annotation only
\usepackage{algorithm}
\usepackage[noend]{algpseudocode}
\usepackage{booktabs}
\usepackage{subcaption}
\usepackage{multicol}
\usepackage[numbers]{natbib}
\usepackage[margin=1in]{geometry}
\usepackage{amsthm}
\newtheorem{definition}{Definition}
\newtheorem{remark}{Remark}
\definecolor{psd}{RGB}{0,0,0}           %0
\definecolor{dd}{RGB}{231,  49, 51}      %1
\definecolor{cdd}{RGB}{75, 139, 191}    %2
\definecolor{psd-dd}{RGB}{95,183, 92}   %3
\definecolor{dd-psd}{RGB}{255, 139, 26} %4
\definecolor{sdd}{RGB}{220, 207, 110}   %5
\definecolor{csdd}{RGB}{175, 103, 62}   %6
\definecolor{dd-star}{RGB}{163, 163, 163}   %8

\definecolor{reg0}{RGB}{51, 143, 188}
\definecolor{reg1}{RGB}{254, 193, 110}
\definecolor{reg2}{RGB}{246, 122, 69}
\definecolor{reg3}{RGB}{217, 67, 79}

\DeclareMathOperator*{\argmin}{arg\!\,min}

\newcommand{\R}{\mathbb{R}} 
\newcommand{\K}{\mathbb{K}}

\newcommand{\Z}{\mathbb{Z}}

\newcommand{\A}{\mathcal{A}}

\newcommand{\psd}{\mathbb{S}}
\newcommand{\cs}{\mathcal{C}}
\newcommand{\ks}{\mathcal{K}}
\newcommand{\gs}{\mathcal{G}}
\newcommand{\es}{\mathcal{E}}
\newcommand{\vs}{\mathcal{V}}

\newcommand{\di}{\mathcal{D}}
\newcommand{\dd}{\mathcal{D}\mathcal{D}}

\newcommand{\sdd}{\mathcal{S}\mathcal{D}\mathcal{D}}
\newcommand{\fw}{\mathcal{F}\mathcal{W}}

\DeclarePairedDelimiter{\abs}{\lvert}{\rvert}
\DeclarePairedDelimiter{\norm}{\lVert}{\rVert}

\DeclarePairedDelimiter{\diag}{\text{diag}(}{)}

\DeclarePairedDelimiterX{\inp}[2]{\langle}{\rangle}{#1, #2}

\newcommand{\tr}{{\mathsf T}}

\newtheorem{theorem}{Theorem}
\newtheorem{proposition}{Proposition}

\newcommand{\new}[1]{\textcolor{black}{#1}}

\title{Decomposed Structured Subsets for Semidefinite and Sum-of-Squares Optimization} %

\author{Jared Miller \thanks{J. Miller and M. Sznaier are with the ECE Department, Northeastern University, Boston, MA 02115. (Emails: miller.jare@northeastern.edu,  msznaier@coe.neu.edu).}
\thanks{J. Miller and M. Sznaier were partially supported by NSF grants  CNS--1646121, CMMI--1638234, IIS--1814631 and ECCS—1808381  and AFOSR grant FA9550-19-1-0005.}, 
Yang Zheng \thanks{Y. Zheng is with the Department of Electrical and Computer Engineering, UC San Diego, La Jolla, CA, 92093, US. (Email: zhengy@eng.ucsd.edu).},  Mario Sznaier\footnotemark[1] \footnotemark[2], Antonis Papachristodoulou% <-this % stops a space
\thanks{A. Papachristodoulou is with the Department of Engineering Science, University of Oxford, Oxford, UK OX1 3PJ. (Email: antonis@eng.ox.ac.uk).}
}
\date{}
\begin{document}
 \maketitle

\begin{abstract}
Semidefinite programs (SDPs) are standard convex problems that are frequently found in control and optimization applications. Interior-point methods can %approximately 
solve SDPs in polynomial time up to arbitrary accuracy, but scale
poorly as the size of matrix variables and the number of constraints increases. 
To improve scalability, SDPs can be approximated with lower and upper bounds through the use of structured subsets (e.g., diagonally-dominant and scaled-diagonally dominant matrices).
Meanwhile, any underlying sparsity or symmetry structure may be leveraged to form an equivalent SDP with smaller positive semidefinite constraints. 
In this paper, we present a notion of \emph{decomposed structured subsets} to approximate an SDP with structured subsets after an equivalent conversion.
The lower/upper bounds found by approximation after conversion 
become tighter than the bounds obtained by approximating the original SDP directly. We apply decomposed structured subsets to semidefinite and sum-of-squares optimization problems with examples of $\mathcal{H}_\infty$ %H-infinity 
norm estimation and constrained polynomial optimization. 
An existing basis pursuit method is adapted into this framework to iteratively refine bounds.
\end{abstract}

% \input{paper_sections_extended/section_0_abstract.tex}

% \end{frontmatter}

% \input{paper_sections_extended/section_1_intro.tex}
 \section{Introduction}
Semidefinite programs (SDPs) are a class of convex optimization problems with a linear objective, affine constraints, and an additional positive semidefinite (PSD) constraint on the decision variable. SDPs include common optimization problems such as Linear Programs (LPs) and Second-order Cone Programs (SOCPs). A more general conic program has a cost $C \in \psd_n$, constraint matrices $A_1,\ldots, A_m \in \mathbb{S}^{n}$, and constraint values $b\in\mathbf{R}^m$. Variables are restricted to a proper cone $K$ and dual cone $K^*$, where $\inp{\cdot}{\cdot}$ denotes the canonical inner product between elements in cones. A conic program has the following primal and dual forms:

% A more general conic program has the following primal and dual forms:
% \begin{subequations}
% 
\noindent\begin{minipage}{.5\linewidth}
  \begin{align}\label{Eq:SDPprimal}
   p^* = & \min_{X} \;  \inp{C}{X} \nonumber\\
    \new{\text{s.t.}} \; &b_i = \inp{A_i}{X} \mid_{i=1}^m  \\
    &X \in K\new{,} \nonumber
  \end{align}
\end{minipage}%
\begin{minipage}{.5\linewidth}
  \begin{align}\label{Eq:SDPdual}
    d^* = &\max_{y, Z} \; \inp{b}{y} \nonumber \\
   \new{\text{s.t.}} \; &C = Z \!+\! \textstyle\sum_{i=1}^m y_i A_i  \\
    & y \in \R^m, \; Z \in K^*\new{.} \nonumber
  \end{align}
\end{minipage}

% \end{subequations}
%
% \noindent\begin{minipage}{.5\linewidth}
% \begin{align}%[t]
%     p^* =\min_{X} \quad & \inp{C}{X} \nonumber \\
%     & \inp{A_i}{X} = b_i\mid_{i = 1}^m,\label{Eq:SDPprimal}\\
%      & X \in K,   \nonumber 
%   \end{align}
% \end{minipage}%
% \begin{minipage}{.5\linewidth}
%   \begin{align}%[t]
%      d^* = \max_{y, Z} \quad & \inp{b}{y} \nonumber\\
%      & Z + \sum_{i=1}^m {y_i A_i}= C, \label{Eq:SDPdual}\\
%      & Z \in K^*, \nonumber
%   \end{align}
% \end{minipage}
The objectives in~\eqref{Eq:SDPprimal} and \eqref{Eq:SDPdual} are related by $p^* \geq d^*$, which is known as weak duality \cite{boyd2004convex}. Strong duality, where $p^*=d^*$, may hold under appropriate constraint qualification conditions (e.g. Slater). \new{SDPs are conic programs with $K = K^* = \psd_+^n$, where $\psd_+^n$ denotes the set of PSD matrices. The dual form \eqref{Eq:SDPdual} of an SDP} is also known as a Linear Matrix Inequality (LMI) \cite{boyd1994linear}. 

Classical interior-point methods (IPMs) can solve an SDP to $\epsilon$-accuracy in polynomial time with complexity $\mathcal{O}(n^2 m^2 + n^3 m \new{+ m^3})$ per iteration \citep{alizadeh1995interior}.
%
%The computational complexity of solving an $n$-dimensional SDP with $m$ equality constraints is polynomial in $n$ and $m$, specifically for Primal-Dual interior point methods (IPM) it scales as $O(n^2 m^2 + n^3 m)$ \citep{alizadeh1995interior}. 
When $m$ is fixed, %based on the problem (after processing redundant constraints),
the speed of IPMs can be greatly improved by reducing the size of PSD cone $\mathbb{S}^n_+$. This motivates a variety of decomposition methods, which exploit problem structures to break up a large PSD constraint into a product of smaller PSD constraints. For example, sparsity in problem data $(C, A_i)$ motivates a notion of chordal decomposition~\citep{agler1988positive,grone1984positive},
and symmetry/common *-algebra structure of $(C, A_i)$ restricts optimization to an invariant subspace~\citep{vallentin2009symmetry}. 
% A survey of decompositions taking advantage of problem structures was documented in~\cite{de2010exploiting}. 

% If the SDP is generally dense, bounds of the true objective can be obtained by inner and outer approximations \citep{ahmadi2017optimization}. 

Structured subset methods restrict~\eqref{Eq:SDPprimal} to simple subsets $K_{\text{inner}} \subset \psd^n_+ \subset K_{\text{outer}}$ \new{to develop inner and outer approximations resulting in optima $p^*_{\text{outer}} \leq p^*_{\text{SDP}} \leq p^*_{\text{inner}}$.} These simple subsets include (scaled-) diagonally dominant (DD or SDD) cones \citep{barker1975cones, boman2005factor}.  \new{Solving~\eqref{Eq:SDPprimal} where $X$ is DD is an LP, and the scenario where $X$ is SDD is an SOCP}. \new{These simplified formulations yielding possibly conservative bounds are often} much faster to solve than the original SDP. Structured subset approximations can be iteratively refined through a change of basis scheme \cite{hall2018optimization}; see \cite{Majumdar2019ASO} for an overview of decomposition methods and structured subsets in solving SDPs. 
We note that polynomial optimization \new{problems} can be approximated by a hierarchy of sum-of-squares (SOS) programs, which can be cast as structured SDPs~\citep{parrilo2000structured}. 
% The resulting SDPs may inherit symmetry and sparsity properties of the original polynomials~\citep{gatermann2004symmetry}. 
The method of structured subsets has also been used to find bounds on polynomial optimization problems when the standard SOS method leads to prohibitively large SDPs; see~\cite{majumdar2014control}.  
Structured subset techniques (DD/SDD matrices) ignore any underlying sparsity and reducible properties in the original SDP. \new{For example, a} diagonally dominant constraint imposes that each diagonal element is greater than the sum of all absolute values on its row/column. Even if the original problem is sparse (comparatively few elements appear in cost or constraints), the problem approximated by \new{standard} DD/SDD constraints will still be dense and may have a slower runtime than the sparse-converted SDP\new{~\citep{vandenberghe2015chordal}}. On the other side, an SDP with sparse/symmetric structure may \new{still} have overly large blocks \new{after conversion}, and these PSD blocks may dominate the computational performance.%~\cite{sun2014decomposition}

This paper presents the notion of \emph{Decomposed Structured Subsets} to find improved lower/upper bounds to semidefinite programs, which exploits problem properties (e.g. sparsity and symmetry) before approximating it with a structured subset. 
Eigenvalue tests of the primal/dual solution can be used to certify if an approximation achieves the true SDP optimum, and containment properties of decomposed structured subsets are presented with their effects on the resultant bounds. 
The cones in the decomposition may be mixed, such that large PSD blocks are approximated with structured subsets while small blocks remain PSD \new{to yield tighter bounds than a uniform cone approximation}.

Some preliminary results were presented at the virtual IFAC 2020 world congress in  \new{\cite{miller2020decomposed}}. This paper additionally explores SDPs with multiple kinds of structure, including SDPs that simultaneously have sparsity and symmetry. 
% Using all possible decompositions (symmetry, then sparsity) before approximating leads to the most accurate bounds. 
Decomposed structured subsets are applied to polynomial optimization problems 
% with multiple methods of leveraging sparsity in sum-of-squares programs. 
\new{through sum-of-squares approximations.}
The analysis of the iterative change of basis algorithm \new{to our framework} is extended. 

The rest of this paper is organized as follows. Section \ref{sec:prelim} introduces preliminaries regarding chordal decomposition and structured subsets. Section~\ref{sec:decomposed} unites these concepts with decomposed structured subsets and performs a containment analysis. Section~\ref{sec:semidefinite} discusses how to apply decomposed structured subsets to semidefinite programs and the change of basis algorithm. This approach is demonstrated through $\mathcal{H}_{\infty}$ norm estimation of networked systems in Section \ref{sec:hinf}. The extension to SOS optimization is covered in Section \ref{sec:polynomial}. 
We conclude this paper in Section~\ref{sect:conclusion}. Appendix \ref{Sec:appendix} contains details of DD/SDD decompositions, and Appendix \ref{sec:symm_appendix} extends the decomposed structured subset framework to problems with algebraic symmetry.

\section{Preliminaries}
\label{sec:prelim}

%This section will review Structured Subsets of the PSD Cone, chordal graphs and decompositions, and sum-of-squares optimizaition.

\subsection{Structured Subsets}

\label{sec:structured}

A basic structured subset of the PSD cone $\mathbb{S}^n_+$ is diagonal nonnegative matrices $\di$.  Two additional subsets are the cones of diagonally dominant (DD) \citep{barker1975cones} and scaled diagonally dominant (SDD) matrices \citep{boman2005factor}:
\begin{align}
    %\di^n &= \{A  \in \psd^n: \exists a_{i} \geq 0 \mid A = \diag{a} \} \\
    \di^n &= \{A  \in \psd^n:  A = \diag{a_1, \ldots, a_n}, a_{i} \geq 0 \}, \nonumber\\
    \dd^n &= \{A \in \psd^n: \; a_{ii} \geq \textstyle\new{\sum_{j \neq i} }{\abs{a_{ij}}}, i = 1, 2, \ldots, n\}, \label{Eq:subsetdef}\\
    \sdd^n &=  \{A \in \psd^n: \; \exists D \in \di^n \; \vert \;  D A D \in \dd^n\}. \nonumber
\end{align}
These subsets satisfy the following containment relation
\begin{equation} \label{eq:ddsddpsd}
    \di^n \subset \dd^n \subset \sdd^n \subset \mathbb{S}^n_+.
\end{equation}

Optimizing conic program \eqref{Eq:SDPprimal} by setting $K$ equal to these cones with a minimization objective will find bounds:
\begin{equation} %\label{eq:ddsddpsd_bounds}
    p_\di \geq p_{\dd} \geq p_{\sdd} \geq p_{\text{SDP}}.
\end{equation}

Solving the conic program over $\di^n$ and $\dd^n$ (\emph{i.e.}, setting $K=\di^n$ or $K=\dd^n$ in~\eqref{Eq:SDPprimal}) is an LP, and over $\sdd^n$ (\emph{i.e.}, setting $K=\sdd^n$ in~\eqref{Eq:SDPprimal}) is an SOCP \citep{ahmadi2017dsos}. As there exist very efficient solvers for LPs and SOCPs, these inner approximations to SDPs can scale to very large-dimension problems.
% Figure \ref{fig:dd-star} shows a PSD-representable feasible set (black), along with inner approximation found by optimizing over the cone $\dd^n$ (red). An outer (gray) LP approximation is also shown, corresponding to the dual cone $(\dd^n)^*$. This example is discussed in Section \ref{Sec:M_feasibility}.  
% The solution quality will necessarily be conservative; see extensive numerical examples in~\citep{ahmadi2017dsos}.

% \begin{figure}[ht]
%     \centering
%     \includegraphics[width=0.4\linewidth]{img/clean_dd_star.pdf}
%     \caption{Cone slices of~\eqref{eq:exampleM} showing ${\color{dd}\dd^n} \subset  \psd_+^n \subset {\color{dd-star} (\dd^n)^*}$: the red region corresponds to solutions of~\eqref{eq:exampleM} over $\dd^n$; the black region corresponds to solutions of~\eqref{eq:exampleM} over $\psd_+^n$; the gray region corresponds to solutions of~\eqref{eq:exampleM} over $(\dd^n)^*$. }
%     \label{fig:dd-star}
% \end{figure}

\label{sec:factor_width}
Factor width matrices also form a structured subset of $\mathbb{S}^n_+$. A matrix $M \in \fw^n_k$ if there exists a rectangular matrix $U$ such that $M = U U^\tr$ where each column of $U$ has cardinality at most $k$ \citep{boman2005factor}. 
%Equivalently:
%\[ \fw^n_k = \left\{ A | \ A = \sum_{I \subset \binom{[n]}{k}} {\varepsilon^n_I(A_I)} \quad  A_I \in \psd_+^\ell \right\} \]
%Notation here is from \citep{ding2018higher}, where $\binom{[d]}{k}$ is the set of rising $k$-length sequences from $1\ldots d$, $ \varepsilon ^n_I(A_I)$ embeds the $k\times k$ matrix $A_I$ into a $d \times d$ matrix as indexed by $I$ with all other values 0, and its left inverse $\tau_I$ extracts out the $k \times k$ submatrix of $A$ as indexed by $I$. 
An intuitive interpretation is that factor width-$k$ matrices are the sum of $k \times k$ PSD matrices that are embedded in $n \times n$ larger matrices. Factor width matrices can be extended to partitions of indices. A block factor-width $k$ matrix \new{given a partition of indices }
is a matrix $M = U U^T$ where each column of $U$ has nonzero elements in at most $k$ sets in the partition \citep{zheng2019block}. In Section~\ref{sec:hinf}, $B_k$ is defined as the set of block factor-width 2 matrices where each block is of size $k$ (up to divisibility).

\subsubsection{Change of Basis}
\label{Sec:ChangeofBasis}
% The notion of change of basis is an iterative method that refines an existing structured subset for approximating SDPs \citep{ahmadi2015sum}. 
The change of basis method is an iterative algorithm that sharpens bounds from structured subsets \citep{ahmadi2015sum}.
%The logic of the change of basis is as follows. 
% Let $\mathcal{P}$ be a property, e.g., the property that a matrix is diagonal or $SDD$. If a matrix $X$ does not satisfy $\mathcal{P}$, there may exist a basis transformation $B X B^\tr $ that does have property $\mathcal{P}$. To be concrete, let $\mathcal{P}$ to denote that a matrix is diagonal. A generically symmetric matrix $X$ is not diagonal, while $B X B^\tr $ satisfies $\mathcal{P}$ if the basis $B$ is the transpose of the eigenvectors of $X$. 
% The Change of Basis algorithm extends this logic to structured subsets. 
% Properties such as a matrix being $DD$ may not be preserved by basis transformation
% As an example, a matrix $X$ may not be $DD$, but $B X B^\tr $ is $DD$ for a suitably chosen $B$.
Given a basis-change matrix $B \in \R^{n \times n}$ and a structured subset cone $K \subset \psd^n_+$, the basis-changed cone is %$K^n(B) = \{X \mid B X B^\tr \in K^n\}$. 
 $K(B) = \{B Q B^\tr \mid Q \in K\}$.
PSD matrices $X \not\in K$ can be made $X \in K(B)$ \new{for some appropriate basis $B$}. To start the iterative refinement process, we first solve a conic optimization problem over a structured subset $K$ such as in \eqref{Eq:SDPprimal}, leading to the iterate $X_0$. 
\new{The Cholesky decomposition $X_0 = L_0L_0^\tr$ can be used to find the next optimal solution $X_1$:}
% Then, we use $X_0 = L_0L_0^\tr$ to formulate a basis-change matrix $B:= L_0$%For example, Hall et. al. chose lower triangular bases $B$ based on Cholesky factors, referred to as $L_0$
% . %If $X_0 \in K$ and $X_0$ has full rank, then $X_0 = L_0 I L_0^T$ and $I \in K(L_0)$. 
% This Cholesky decomposition can be used  to find the next optimal solution $X_1$:
%
%
% Given a matrix $L$ and a cone $K$, a basis-changed cone is $K(L) = \{X \mid L X L^\tr \in K\}$.  After finding an optimum $X_0$ of the conic optimization \eqref{Eq:SDPprimal}, form a matrix factorization $X_0 = L_0 L_0^\tr$, and then solve the modified problem:
%
\begin{equation}
    \label{Eq:ChangePrimal}
    \begin{aligned} 
        X_1 = \argmin_{X} \quad & \inp{C}{X} \\
    \new{\text{s.t.}} \quad & \inp{A_i}{X} = b_i, i = 1, \ldots, m, \\
    & X \in K(L_0).
    \end{aligned}
\end{equation}
Use $X_1 = L_1 L_1^T$, and solve the same problem over $K^n(L_1)$ to find optimal point $X_2$. The cost $\inp{C}{X}$ of iterate $t$ is upper bounded by the cost at iterate $t-1$,
% \new{since the optimal solution $X_{t-1}$ naturally belongs to $K(L_{t-1})$ that is the feasible region of the next iteration.}
\new{because both $X_{t-1}$ and $X_t$ are members of the feasible set $K(L_{t-1})$ at iteration $t$.}
% since $I \in K(L_{t-1})$ is a feasible point. %The volume of each cone $K(L_t)$ decreases between iterations, and 
If the structured subset chosen $K^n = \dd^n$, then the SDP is approximated by an iterative sequence of LPs. 
We note that \new{this} procedure might not converge to the true SDP optimum. 
An analogous process can occur on the dual side to find an increasing sequence of lower bounds to the true SDP cost; see~\cite{hall2018optimization} for details. 
% A column generation process to refine structured subset approximations also exists .

\subsection{Chordal Decomposition}

\label{sec:chordal}

An SDP is sparse if only a few entries of $X \in \psd_+^n$ are involved in the cost and constraints. For example, if $C_{jk} = (A_{i})_{jk} = 0, \  \forall i = 1, \ldots, m$, the value\new{s} of $X_{jk}$ and $X_{kj}$ \new{are} simply present to ensure that $X \succeq 0$. The aggregate sparsity pattern of $(C, A_i)$ can be encoded by a graph $\gs(\vs, \es)$, where there is an edge between vertices $i$ and $j$ if any of $C, A_1,\ldots ,A_m$ is nonzero at indices $(i, j)$. We now introduce a few graph-theoretic notions. A chord in a graph is an edge between two non-consecutive vertices in a cycle, and a graph is chordal if every cycle of length 4 or more has a chord \citep{vandenberghe2015chordal}. Non-chordal graphs can be chordal-extended by adding edges, and heuristics exist to approximate minimum fill-in \citep{yannakakis1981computing}. 
A clique $\cs$ is a set of vertices that forms a complete graph: $\forall v_i, v_j \in \cs, (v_i, v_j) \in \es$. Maximal cliques are cliques that are not contained in another clique. The cardinality of a maximal clique is denoted as $|\mathcal{C}|$.  

Following notation from \new{\cite{kakimura2010direct}}, the set of sparse symmetric matrices with pattern $\gs$ forms a cone  
\new{$\psd^n(\es, 0) = \{X \in \psd^n \mid \ X_{ij} = 0, \ \forall i \neq j, \  (i, j) \not\in \es\}$}.
The sparse PSD cone defined by $\es$ is \new{$ \psd^n_+(\es, 0) = \psd^n(\es, 0) \cap \psd_+^n$.} The dual cone $[\psd^n_+(\es, 0)]^* = \psd^n_+(\es, ?)$ is the set of sparse symmetric matrices that admit a PSD completion\new{.}

\new{For a vector $x \in \R^n$ and a clique $\mathcal{C} \subseteq \vs$, there exists a vector $x_{\cs} \in \R^{\abs{\mathcal{C}}}$ that selects the entries of $x$ with indices $\mathcal{C}$.
}
Let $E_{\mathcal{C}} \in \mathbb{R}^{|\mathcal{C}| \times n}$ be $0/1$ entry selector matrices 
% that index out entries in clique $\cs_k$, 
such that $x_{\cs} = E_{\cs} x, \new{\ \forall x \in \R^n}$. The cones $\psd^n_+(\es, ?)$, and $\psd^n_+(\es, 0)$ have a decomposable structure if $\gs$ is chordal:

% Given a graph $\gs(\vs, \es)$, let $\es^*$ be an edge set $\es$ augmented with self-loops. The cone of sparse symmetric matrices $$\psd^n(\es, 0) = \{X \in \psd^n \mid \ X_{ij} = 0, \ \forall (i, j) \not\in \es^*\},$$ and its subcone of sparse PSD symmetric matrices is 
% $$
%     \psd^n_+(\es, 0) = \psd(\es, 0) \cap \psd_+^n.
% $$ 
% The dual space $\psd^n_+(\es, ?)=[\psd^n_+(\es, 0)]^*$ are matrices with entries on $\es^*$ that can be completed into PSD matrices. Let $E_{\mathcal{C}_k} \in \mathbb{R}^{|\mathcal{C}_k| \times n}$ be $0/1$ entry selector matrices that index out entries in clique $\cs_k$. % $\cs_k$ ($x_k = E_k x$). 
% For sparse PSD matrices, we have the following two decomposition results.

\begin{theorem} [Grone's Theorem \citep{grone1984positive}] \label{T:GroneTheorem}
     Let $\mathcal{G}(\mathcal{V},\mathcal{E})$ be a chordal graph with a set of maximal cliques $\{\mathcal{C}_1,\mathcal{C}_2, \ldots, \mathcal{C}_p\}$. Then, $X\in\mathbb{S}^n_+(\mathcal{E},?)$ if and only if
     $$ X_k = E_{\mathcal{C}_k} X E_{\mathcal{C}_k}^\tr \in \mathbb{S}^{\vert \mathcal{C}_k \vert}_+,
    \qquad k=1,\,\ldots,\,p.$$
    \end{theorem}

\begin{theorem} [Agler's Theorem {\citep{agler1988positive}}] \label{T:AglerTheorem}
     Let $\mathcal{G}(\mathcal{V},\mathcal{E})$ be a chordal graph with a set of maximal cliques $\{\mathcal{C}_1,\mathcal{C}_2, \ldots, \mathcal{C}_p\}$. Then, $Z\in\mathbb{S}^n_+(\mathcal{E},0)$ if and only if there exist $Z_k \in \mathbb{S}^{\vert \mathcal{C}_k \vert}_+, k=1,\,\ldots,\,p,$ such that
     \begin{align*}
    Z &= \sum_{k=1}^p{E_{\mathcal{C}_k}^\tr Z_k E_{\mathcal{C}_k}}.
\end{align*}
    \end{theorem}

Theorem~\ref{T:GroneTheorem} breaks up a large sparse PSD constraint $X \in \mathbb{S}^n_+(\mathcal{E},?)$ into a series of smaller coupled PSD constraints $X_k \succeq  0, k =1, \ldots, p$. This result can be applied to primal SDPs with a chordal sparsity pattern $\mathcal{E}$, \emph{i.e.}, problem~\eqref{Eq:SDPprimal} with $K = \mathbb{S}^n_{+}(\mathcal{E},?)$ can be decomposed as % () %More formally, Grone's theorem forms an equivalence between $\psd_+(\es, ?)$ and a fiber product of cones $\psd_+$. This fiber product is conducted with respect to restrictions mappings: cliques must be consistent across their overlaps.
%There is no coupling between clique matrices $Z_k$, but $X_k$ require equality constraints to ensure consistency between overlapping cliques. The chordalized SDP ($X \in \psd^n_+(\es, ?)$ in Primal form is:
%
\begin{equation} \label{Eq:PrimalSDP}
    \begin{aligned}
    \min_{X} \quad & \langle C,X \rangle \\
    \text{subject to} \quad & \langle A_i,X \rangle = b_i, i = 1, \ldots, m, \\
    & E_{\cs_k} X E_{\cs_k}^\tr  \in \mathbb{S}^{|\mathcal{C}_k|}_+, k = 1, \ldots, p.
    \end{aligned} 
\end{equation}
Analogous results can be obtained for sparse dual SDPs, with a characterization of $Z \in \psd^n_+(\es, 0)$ using Theorem~\ref{T:AglerTheorem}. % \cite{agler1988positive}. 
These decomposed SDPs can be solved using first order methods to low accuracy via variable splits $E_{\cs_k} X E_{\cs_k}^\tr = X_k$ (see~\citep{zheng2019chordal} for details), but interior point methods may suffer from the increase of the equality constraints introduced by the decomposition \citep{vandenberghe2015chordal}. Conversion utilities such as SparseCoLO \citep{fujisawa2009user} internally perform domain- and range-space decompositions to exploit chordal sparse structures. An algorithm to trade off between PSD block sizes and added equality constraints was discussed in \cite{garstka2020clique}.

%These conversions may end up with a graph $\gs_F(\vs_F, \es_F) \supset \gs(\vs, \es)$ with more edges than the chordal decomposition (fill-in), and therefore may contain larger cliques.

% \input{paper_sections_extended/section_3_decomposed_structure.tex}

\section{Decomposed Structured Subsets}

\label{sec:decomposed}

This section combines decomposition methods and structured subsets into \emph{decomposed structured subsets}.

\begin{figure}[ht]
    \centering
    \includegraphics[width=0.6\linewidth]{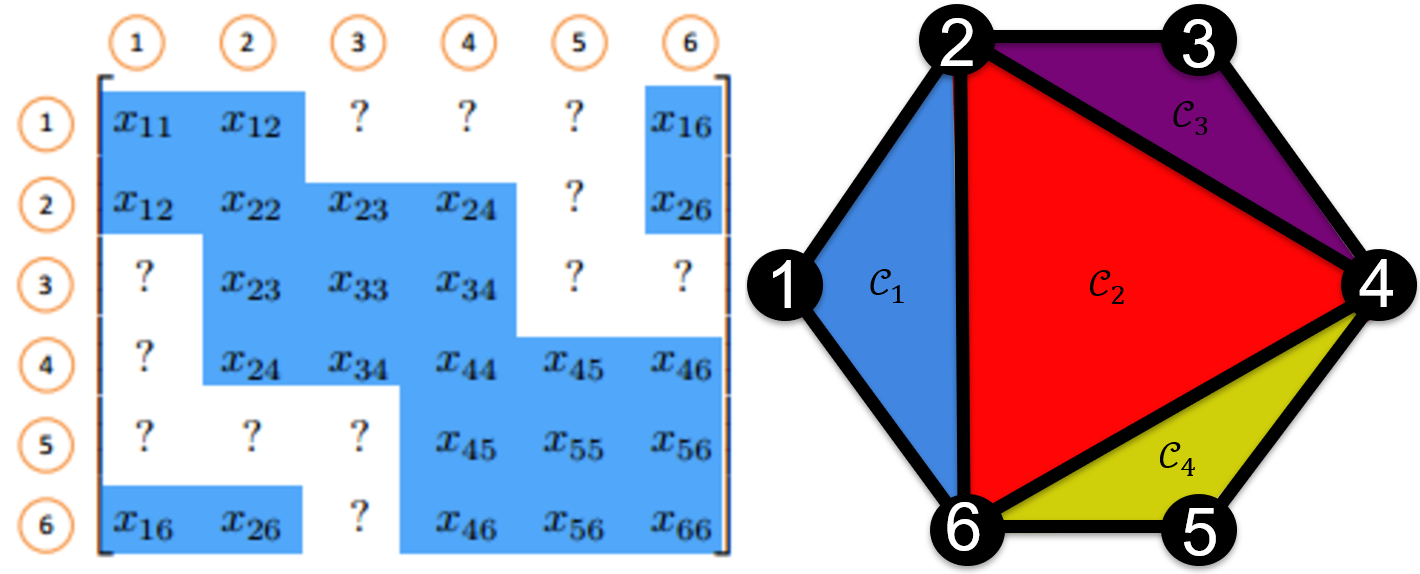}
    \caption{Left: A $6 \times 6$ PSD completable cone  $\mathbb{S}^n_+(\mathcal{E},?)$. 
    %SDP with aggregate sparsity pattern (left), 
    Right: corresponding chordal graph $\gs(\vs,\es)$ with maximal cliques $\{\cs_k\}_{k=1}^4$.}
    \label{fig:grone_hex}
\end{figure}

As an example, consider problem~\eqref{Eq:SDPprimal} with $K = \psd^6_+(\es, ?)$ where $\es$ is the sparsity pattern shown in Figure~\ref{fig:grone_hex}. Theorem~\ref{T:GroneTheorem} poses an optimization problem over the cliques $\{X_k \in \psd_+^3\}_{k=1}^4$. %, while a standard interior point method would take in a variable $X \in \psd_+^6$ where only variables in $\es$ matter in optimization (\emph{e.g.} $x_{13}$ does not appear in constraints nor cost, so it has a `?' entry). 
Now consider a structured subset restriction.
%The effect of restriction to structured subset can be illustrated by examining the variable $x_{11} \in \cs_1$. 
If we require $X = [x_{ij}] \in \dd^6$, this constraint requires $x_{11} \geq \sum_{i=2}^6 {|x_{1i}|}$. Instead, if we consider a decomposition and impose structured subset restriction on the cliques, \emph{e.g.}. $X_1 \in \dd^3$, then it requires $x_{11} \geq \abs{x_{12}} + \abs{x_{16}}$, which is less restrictive than competing against all variables in the same row/column. \emph{Decomposed Structured Subsets} arise from performing decompositions before applying structured subsets, and are presented in detail in this section.
Figure~\ref{fig:grone_hex} shows a chordal graph and its maximal cliques $\mathcal{C}_k, k =1, \ldots, 4$.

\subsection{Definition of decomposed structured subsets}
%Given a chordal graph $\mathcal{G}(\mathcal{V},\mathcal{E})$ with maximal cliques $\mathcal{C}_1, \ldots, \mathcal{C}_p$,

A clique edge cover of a graph $\gs(\vs, \es)$ is a set of subsets $\{\cs_k\}_{k=1}^p$ such that every max. clique of $\gs$ is contained in at least one \new{clique} $\cs_k$. Clique edge covers allow for clique \new{extensions and merges for possibly non-chordal graphs}.

\begin{definition}
\new{We define sparse DD and SDD matrices as:}
% We define z;arse DD and SDD matrices 
$$
\begin{aligned}
    \dd^n(\mathcal{E},0) &=  \mathbb{S}^n(\mathcal{E},0) \new{\ \cap \ \dd^n}, \\
    \sdd^n(\mathcal{E},0) &=  \mathbb{S}^n(\mathcal{E},0)\new{\ \cap \ \sdd^n}.
\end{aligned}
$$
\end{definition}
These sparse matrices obey the containment:
$$
    \dd^n(\mathcal{E},0)  \subset \sdd^n(\mathcal{E},0) \subset \mathbb{S}^n_+(\mathcal{E},0). 
$$
The following result is proven in Appendix \ref{Sec:appendix}:
\begin{proposition} \label{prop:decomposition}
    Let $\mathcal{G}(\mathcal{V},\mathcal{E})$ be a graph with a clique edge cover $\{\mathcal{C}_1,\mathcal{C}_2, \ldots, \mathcal{C}_p\}$. Then,
    \begin{enumerate}
        \item $Z\in \dd^n(\mathcal{E},0)$ iff
     $$
    Z = \textstyle \textstyle \sum_{k=1}^p{E_{\mathcal{C}_k}^\tr Z_k E_{\mathcal{C}_k}}, \quad 
Z_k \in \dd^{\vert \mathcal{C}_k \vert},
    \; k=1,\,\ldots,\,p.
$$
\item 
$Z\in \sdd^n(\mathcal{E},0)$ iff
     $$
    Z = \textstyle \textstyle \sum_{k=1}^p{E_{\mathcal{C}_k}^\tr Z_k E_{\mathcal{C}_k}}, \quad 
Z_k \in \sdd^{\vert \mathcal{C}_k \vert},
    \; k=1,\,\ldots,\,p.
$$
    \end{enumerate}
\end{proposition}
% The proof is provided in Appendix \ref{Sec:appendix}. These results hold for an arbitrary clique edge cover that covers all maximal cliques.

Motivated by Theorems~\ref{T:GroneTheorem} and~\ref{T:AglerTheorem}, and Proposition~\ref{prop:decomposition}, we let $\es$ be a sparsity pattern,  $\ks = \{K_k\}_{k=1}^p$ be a set of cones corresponding to a clique edge cover $\mathcal{C}_1, \ldots, \mathcal{C}_p$, where each individual cone $K_k$ is some structured subset in $\mathbb{S}^{|\mathcal{C}_k|}$.
%of size $\abs{\cs_k}$. 
We  define %two :
two \emph{decomposed structured subsets}:
\begin{equation} \label{eq:decompositionsubsets}
    \begin{aligned}
        \mathcal{K}(\mathcal{E},0) &:= \Big\{Z \in \mathbb{S}^n \mid Z = \sum_{k=1}^p{E_{\mathcal{C}_k}^\tr Z_k E_{\mathcal{C}_k}}, \Big.\;  \\
\Big. &\qquad \qquad \qquad \qquad Z_{_k} \in K_k,
    \; k=1,\,\ldots,\,p\Big\},
        \\
        \mathcal{K}(\mathcal{E},?) &:= \Big\{X \in \mathbb{S}^n \mid  E_{\mathcal{C}_k} X E_{\mathcal{C}_k}^\tr \in K_k,
    \! k=1,\ldots,p\Big\}.
    \end{aligned}
\end{equation}

The decomposed structured subset $\mathcal{K}(\mathcal{E},0)$ generalizes $\dd^n(\mathcal{E},0)$, $\sdd^n(\mathcal{E},0)$, and $\mathbb{S}^n_+(\mathcal{E},0)$ in the following sense:
    \begin{align*}
       \mathcal{K}(\mathcal{E},0) &= \dd^n(\mathcal{E},0), & &   \text{if}\;K_k = \dd^{|\mathcal{C}_k|}, k = 1, \ldots, p. \\
       \mathcal{K}(\mathcal{E},0) &= \sdd^n(\mathcal{E},0), & &  \text{if}\;K_k = \sdd^{|\mathcal{C}_k|}, k = 1, \ldots, p.\\
       \intertext{If $\es$ is chordal, the additional results hold for PSD cones:}
           \mathcal{K}(\mathcal{E},0) &= \mathbb{S}^n_+(\mathcal{E},0), \; & &\text{if}\;K_k = \mathbb{S}^{|\mathcal{C}_k|}_+, k = 1, \ldots, p, \\
           \mathcal{K}(\mathcal{E},?) &=  \mathbb{S}^n_+(\mathcal{E},?), \; & &\text{if}\;K_k = \mathbb{S}^{|\mathcal{C}_k|}_+, k = 1, \ldots, p.
       \end{align*}

\subsection{Containment Analysis}

% The notion of decomposition structured subsets~\eqref{eq:decompositionsubsets} gives more freedom to choose the individual cones $K_k$. We give a detailed analysis below. 
% %
% In Theorems~\ref{T:GroneTheorem} and~\ref{T:AglerTheorem} and Proposition~\ref{prop:decomposition}, the cones corresponding to each clique %cliques 
% are of the same type, \emph{i.e.}, $K_k$ are all either  $\psd_+^{|\mathcal{C}_k|}$ or $\dd^{|\mathcal{C}_k|}$. Additional freedom can be gained by allowing cliques to reside in different cones. 

% Given a graph with cliques $\mathcal{C}_1, \ldots, \mathcal{C}_p$, we consider two sets of cones $\ks = \{K_k\}_{k=1}^p$ and $\tilde{\ks} = \{\tilde{K}_k\}_{k=1}^p$, where $K_k$ or $\tilde{K}_k$ is a cone in $\mathbb{S}^{|\mathcal{C}_k|}$. We define the partial ordering $\subseteq$ on decomposed structured subsets:
% $$
% \ks \subseteq \tilde{\ks} \quad   \text{iff}  \quad  K_k \subseteq \tilde{K}_k  \quad \forall\; k=1 \ldots p.
% $$
% %
% Then, we have the following proposition.
% \begin{proposition} \label{prop:containment}
% Given a graph with cliques $\mathcal{C}_1, \ldots, \mathcal{C}_p$, and two sets of cones $\ks = \{K_k\}^p_{k=1}, \tilde{\ks} = \{\tilde{K}_k\}^p_{k=1}$, if $\ks \subseteq \tilde{\ks}$, then we have
% $$
%     \mathcal{K}(\mathcal{E},0) \subseteq \tilde{\mathcal{K}}(\mathcal{E},0) \quad \text{and} \quad \mathcal{K}(\mathcal{E},?) \subseteq \tilde{\mathcal{K}}(\mathcal{E},?)
% $$
% \end{proposition}

\begin{definition} \label{def:containment_coneset}
For a graph $\gs(\vs, \es)$ with clique cover $\mathcal{C}_1, \ldots, \mathcal{C}_p$, let $\ks = \{K_k\}_{k=1}^p$ and $\tilde{\ks} = \{\tilde{K}_k\}_{k=1}^p$ be two sets of clique-cones. The partial ordering $\subseteq$ is defined as:
%
% Given a graph with cliques $\mathcal{C}_1, \ldots, \mathcal{C}_p$, we consider two sets of cones $\ks = \{K_k\}_{k=1}^p$ and $\tilde{\ks} = \{\tilde{K}_k\}_{k=1}^p$, where $K_k$ or $\tilde{K}_k$ is a cone in $\mathbb{S}^{|\mathcal{C}_k|}$. We define the partial ordering $\subseteq$ on decomposed structured subsets:
$$
\ks \subseteq \tilde{\ks} \quad   \text{iff}  \quad  K_k \subseteq \tilde{K}_k  \quad \forall\; k=1 \ldots p.
$$
\end{definition}

%
% Then, we have the following proposition.
\begin{remark} \label{prop:containment}
% Given a graph with cliques $\mathcal{C}_1, \ldots, \mathcal{C}_p$, and two sets of cones $\ks = \{K_k\}^p_{k=1}, \tilde{\ks} = \{\tilde{K}_k\}^p_{k=1}$, if $\ks \subseteq \tilde{\ks}$, then we have
If clique-cones $\ks \subset \tilde{\ks}$ for a clique edge cover $\mathcal{C}_1, \ldots, \mathcal{C}_p$ of a graph $\gs(\vs, \es)$, then by definition, we have 
$$
    \mathcal{K}(\mathcal{E},0) \subseteq \tilde{\mathcal{K}}(\mathcal{E},0) \quad \text{and} \quad \mathcal{K}(\mathcal{E},?) \subseteq \tilde{\mathcal{K}}(\mathcal{E},?)\new{.}
$$
% This is true from the definition of the partial ordering.
\end{remark}
% The result is true by definition. 
%This remark is true from the definition of the partial ordering and containment of cones. 
%
%
\new{The relationship above is simple yet has useful implications in semidefinite optimization. In particular, t}he cone $\dd^n(\es, 0)$ has the same cone on each clique $\cs_k$ with $X_k \in \dd^{\abs{\cs_k}}$. Mixing cones with $\ks: K_k \supseteq \dd^{\abs{\cs_k}}, \forall k = 1 \ldots p$ where not all $K_k = \dd^{\abs{\cs_k}}$ will form a cone $\ks(\es, 0) \supset \dd(\es, 0)$. Mixing cones therefore results in cones closer to $\psd^n_+(\es, 0)$. Similar statements hold for $\psd_+^n(\es,?)$. \new{This allows us to get better lower and upper bounds for problem~\eqref{Eq:SDPprimal}.}  
% In the context of optimization, DD/SDD constraints offer scalable computation while PSD constraints are close (or exactly meet when the underlying graph is chordal) the true feasible region. Proposition~\ref{prop:containment} suggests some flexibility of choosing the individual cones $K_k$. 
\label{Sec:M_feasibility} 
As an example, consider the following matrix parameterized by  $(a,b)$:
%
%\[
\begin{equation} \label{eq:exampleM}
    M(a,b) = \begin{bmatrix}
         1 & \frac{1}{2}+a & ? & ?\\
         \frac{1}{2}+a & 2 & {-2a} & a+b\\
         ? & {-2a} & 5 & \frac{b}{2}\\
         ? & a+b & \frac{b}{2}  & 2
     \end{bmatrix},    
\end{equation}
where $?$ denotes unspecified entries. The sparsity pattern of 
     %\]
%
$M(a,b)$ has two maximal cliques: $\{2,3,4\}$ and $\{1,2\}$. By Theorem~\ref{T:GroneTheorem}, $M(a,b) \in \psd_+^4(\es, ?)$ if:
\begin{align*}
M_1(a,b) &= \begin{bmatrix}
          2 & {-2a} & a+b\\
          {-2a} & 5 & b/2\\
          a+b & b/2  & 2
     \end{bmatrix} \succeq 0  \\  
     M_2(a,b) &= \begin{bmatrix}
         1 & 1/2+a \\
         1/2+a & 2
     \end{bmatrix}\succeq 0. 
     \end{align*}
Given a cone set $\ks=\{\ks_1, \ks_2\}$, we define the feasibility set for $M \in \ks(\mathcal{E},?)$ as 
$$\{(a,b) \mid M_1(a,b) \in \ks_1 \ \text{and} \  M_2(a,b) \in \ks_2\}.$$

\begin{figure}[ht]
    \centering
    \begin{subfigure}[t]{0.35\linewidth}
        \includegraphics[width=\linewidth]{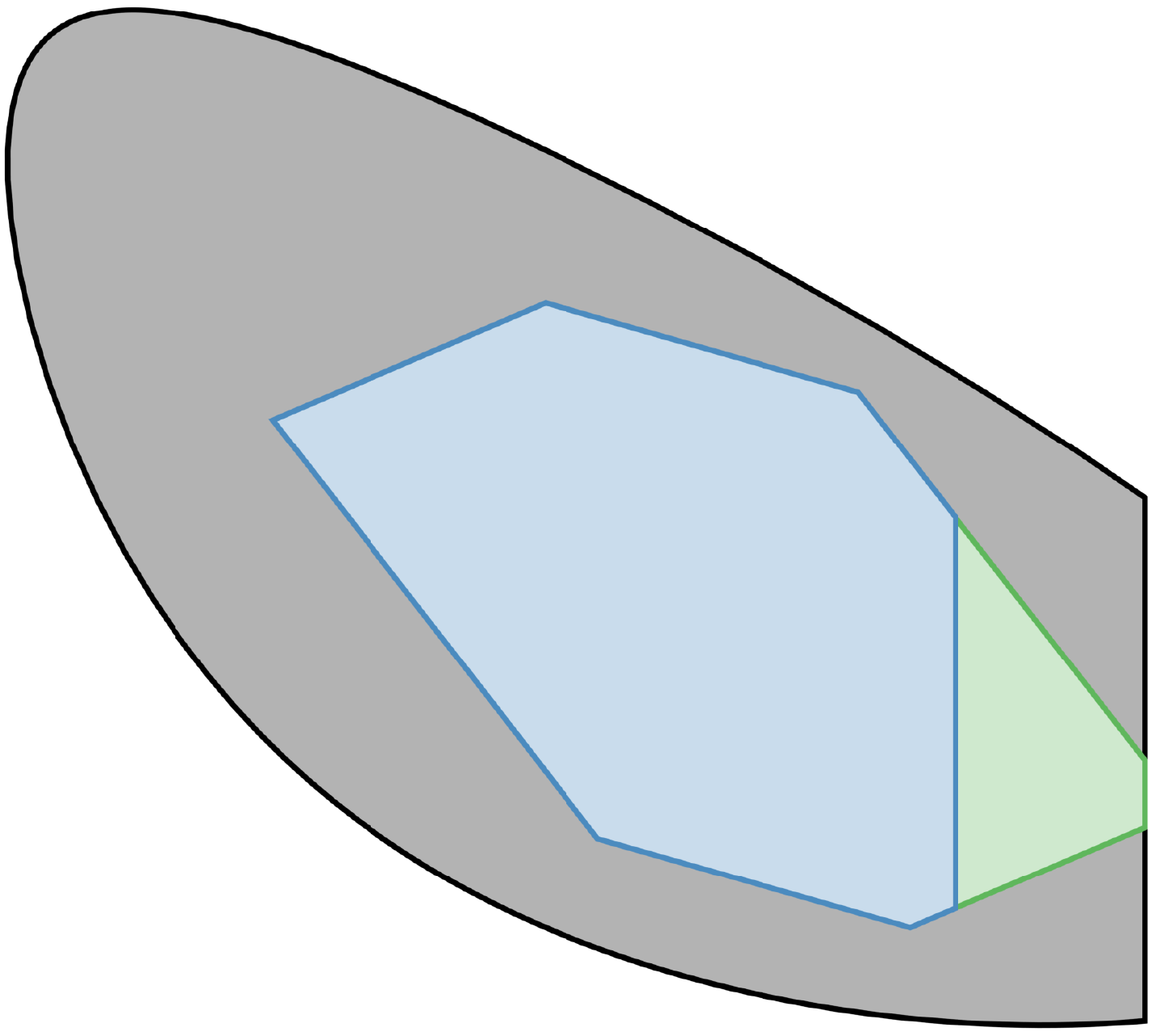}
        %\caption*{${\color{psd-dd} \ks^1 = \{\psd_+^3, \dd^2\}}$}
       % \caption*{${\color{psd-dd} \ks^1 = \{\dd^3,\psd_+^2\}}$}
        \label{fig:mix_dd_psd}
    \end{subfigure}
     \qquad
    \begin{subfigure}[t]{0.35\linewidth}
        \includegraphics[width=\linewidth]{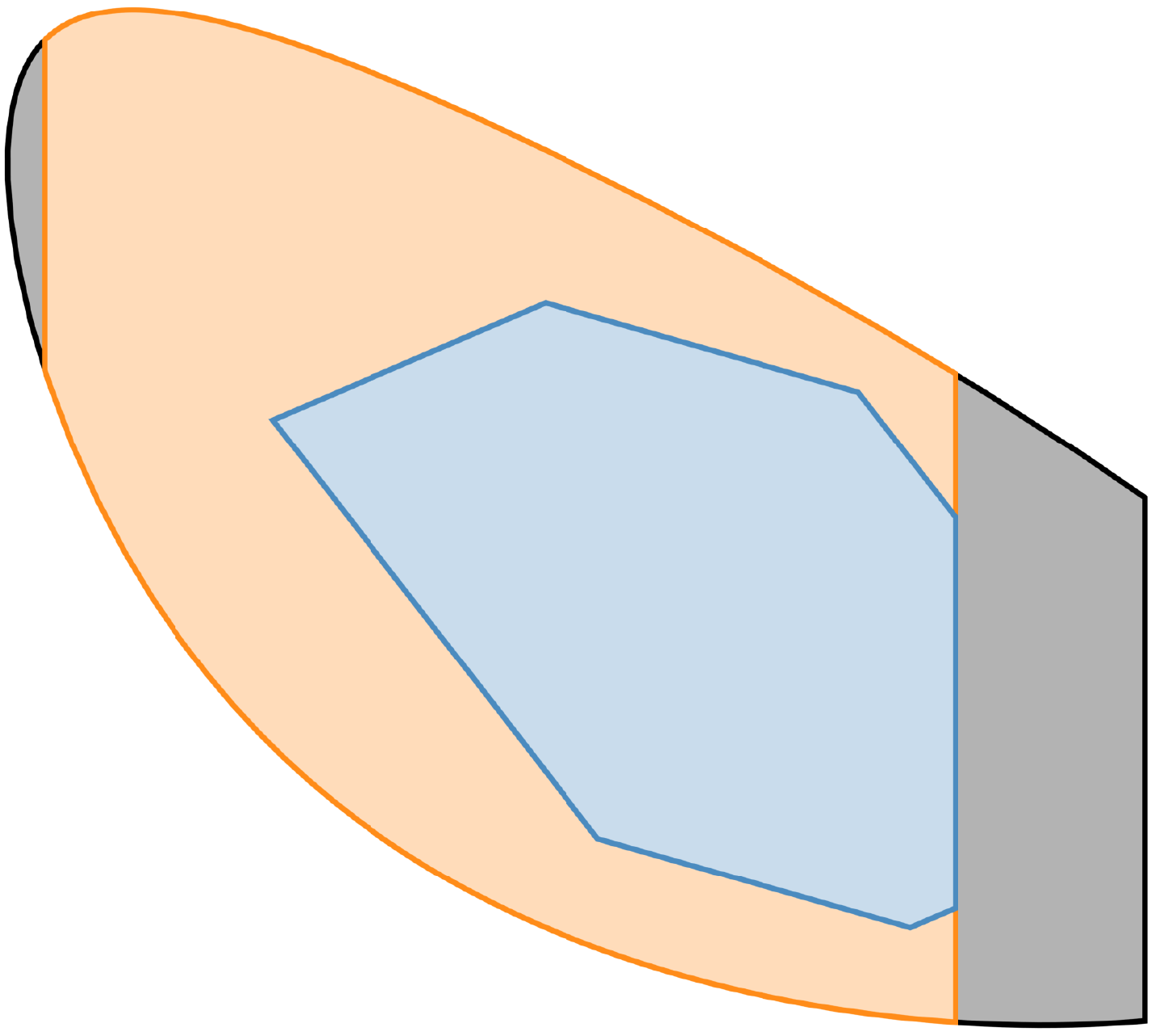}
        %\caption*{${\color{dd-psd} \ks^2=\{\dd^3, \psd_+^2\}}$}
      %  \caption*{${\color{dd-psd} \ks^2=\{\psd_+^3, \dd^2\}}$}
        \label{fig:mix_psd_dd}
    \end{subfigure}
    \caption{Mixing cones broadens feasibility regions for $M \in \ks(\es, ?)$ in~\eqref{eq:exampleM}. Left: ${\color{psd-dd} \ks^1 = \{\dd^3,\psd_+^2\}}$; Right: ${\color{dd-psd} \ks^2=\{\psd_+^3, \dd^2\}}$.}
 \label{fig:cone_mix}
\end{figure}

The black sets in Figure \ref{fig:cone_mix} are subsets of the $(a, b)$-plane where the matrix $M(a,b)$ has a PSD-completion ($M(a,b) \in \psd_+(\es, ?)$). The blue feasibility sets are regions where $M_1(a,b) \in \dd^4$ and $M_2(a,b) \in \dd^2$.
 As expected, the blue set is contained within the black set since $\dd^4(\es, ?) \subset \psd_+^4(\es, ?)$. %$\dd(\es, ?) \subset \sdd(\es, ?)$. 
The orange set in the right panel has %$M_1(a,b) \in \psd^3$ 
$M_1(a,b) \in \psd^3_+$ 
and $M_2(a,b) \in \dd^2$. Note how the orange set includes the blue set (all $\dd$) and expands to nearly fill the left side of the black set (all $\psd_+$). The green set in the left panel has $M_2(a,b) \in \psd_+^2$ instead, which expands the $\dd$ blue set with a small rightward bump.

\begin{definition}
A sparse matrix $X \in \new{\psd}^n(\es, 0)$ has a \textbf{$K$-completion} for a structured subset $K \subseteq \psd_+^n$ if there exists an $\bar{X} \in K$ such that $X_{ij} = \bar{X}_{ij}, \ \forall(i,j) \in \es$. 
% \begin{equation}
% \end{equation}
\end{definition}

\begin{figure}[ht]
    \centering
    \begin{subfigure}[t]{0.35\linewidth}
        \includegraphics[width=\linewidth]{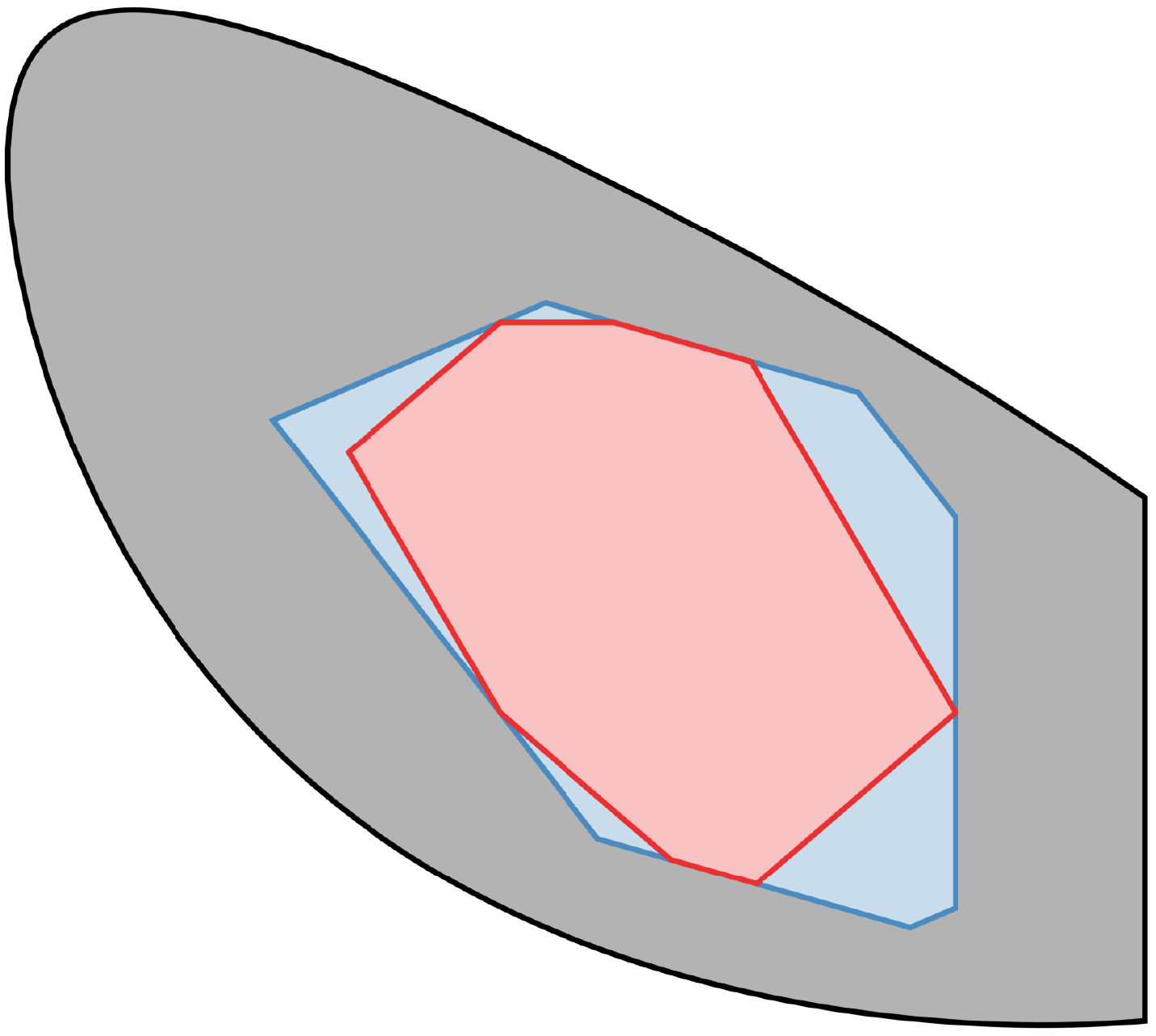}
      %  \caption*{${\color{dd} \dd} \subset {\color{cdd} \dd(\es,  ?)} \subset \psd_+$}
        \label{fig:dd_cdd}
    \end{subfigure} \qquad
        \begin{subfigure}[t]{0.35\linewidth}
        \includegraphics[width=\linewidth]{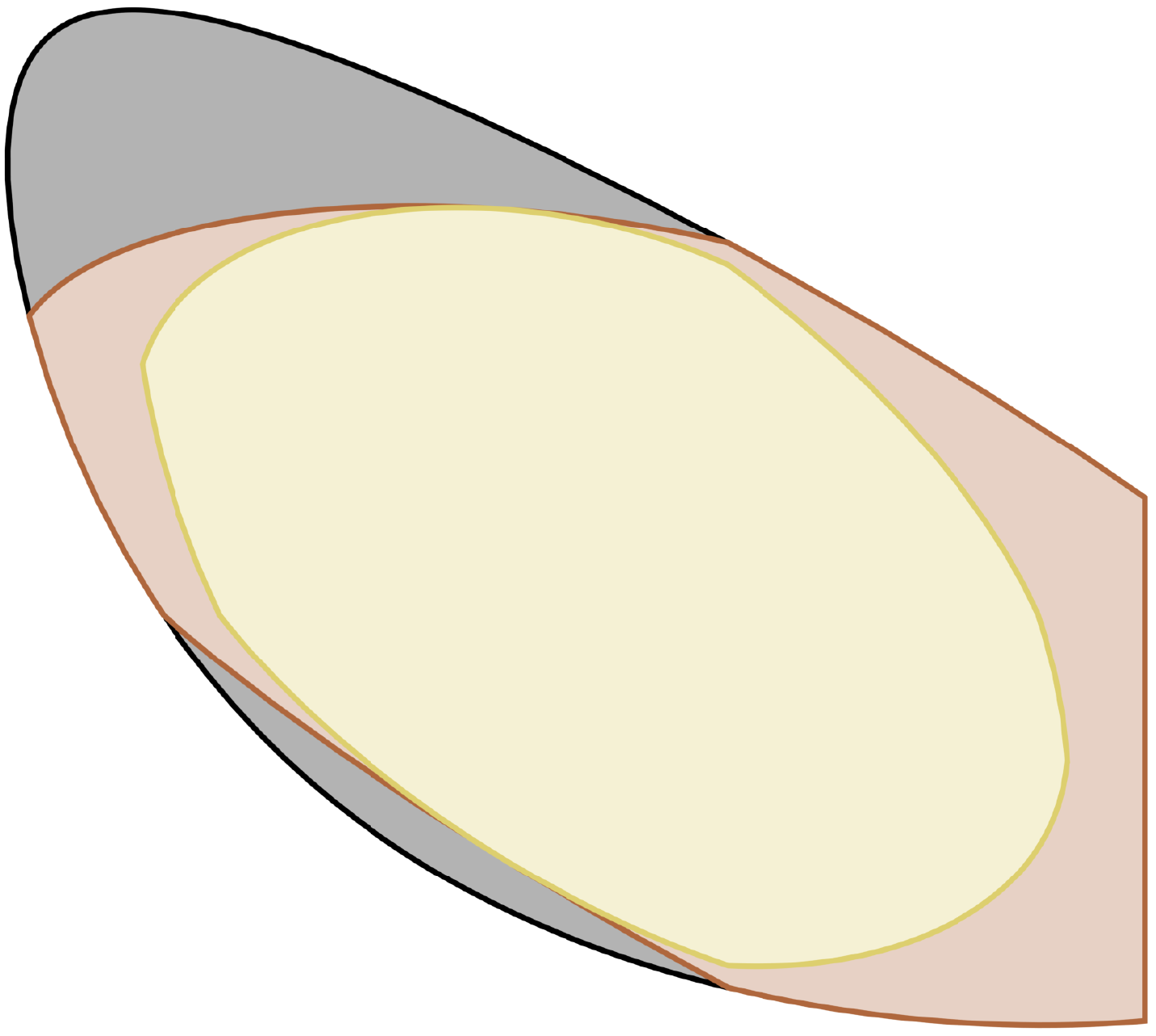}
    %    \caption*{${\color{sdd} \sdd} \subset {\color{csdd} \sdd(\es, ?)} \subset \psd_+$}
        \label{fig:sdd_csdd}
    \end{subfigure}
    \caption{ \label{fig:cone_same}    Regions of the $(a,b)$ plane for which the matrix $M(a,b)$ from Eq. \eqref{eq:exampleM} is ${\color{dd} \dd^4}$-completable, ${\color{cdd} \dd^4(\es,  ?)}$, or $\psd^4_+(\es,?)$ (Left), and ${\color{sdd} \sdd^4}$-completable, ${\color{csdd} \sdd^4(\es, ?)}$, or $\psd_+^4(\es,?)$ (Right).  
    The containment of the plotted regions reflect the cone inclusions
    ${\color{dd} \dd^4}$-completable $ \subset {\color{cdd} \dd^4(\es,  ?)} \subset \psd^4_+(\es,?)$ (Left); ${\color{sdd} \sdd^4}$-completable $ \subset {\color{csdd} \sdd^4(\es, ?)} \subset \psd_+^4(\es,?)$ (Right). }
\end{figure}

%\begin{remark}
\new{For a structured subset $K \subseteq \psd_+^n$ and a sparsity pattern $\gs(\vs,\es)$, the set of $K$-completable matrices with pattern $\es$ is contained within $K(\es, ?)$.}
Figure \ref{fig:cone_same} illustrates and compares feasibility sets for $M(a,b) \in K(\es, ?)$ (cliques of $M(a,b)$ in $K$) and when $M(a,b)$ has a $K$-completion. The blue $\dd^4(\es, ?)$ and black $\psd_+(\es, ?)$ feasibility set are the same in Figure \ref{fig:cone_same} as in \ref{fig:cone_mix}.
The left panel additionally shows the feasible regions where $M(a,b)$ is $\dd^4$-completable (red). The  right plot echoes the left plot, where imposing that $M_1(a,b)\in \sdd, M_2(a,b) \in \sdd$ yields a broader feasibility set than requiring that $M(a,b)$ has an $\sdd$-completion.

Decomposed structured subsets can be posed over dual cones that are larger than the PSD cone.
\begin{proposition}
Let $\mathcal{K} = \{K_k\}_{k=1}^p $ be a set of cones with dual $\mathcal{K}^* = \{K^*_k\}^p_{k=1}$, where $K_k$ is a structured subset in $\psd^{|\mathcal{C}_k|}_+$ (\emph{i.e.}, $\dd^{|\mathcal{C}_k|}$, $\sdd^{|\mathcal{C}|_k}$,  or $\psd^{|\mathcal{C}_k|}_+$), and $\mathcal{C}_k$ is a set of maximal cliques for $\mathcal{E}$. Then, we have 
\begin{align}
    [\mathcal{K}(\mathcal{E},0)]^* &= \mathcal{K}^*(\mathcal{E},?), \label{eq:dual_construction}\\ 
    \qquad [\mathcal{K}(\mathcal{E},?)]^* &\subseteq \mathcal{K}^*(\mathcal{E},0). \label{eq:dual_completion}
\end{align}
\end{proposition}

For the equivalence in Equation \eqref{eq:dual_construction}:
\begin{proof}
Recall the definition of $\mathcal{K}(\mathcal{E},0)$ and $\mathcal{K}(\mathcal{E},?)$ in~\eqref{eq:decompositionsubsets}. We now verify that 
$$
\begin{aligned}
  &[\mathcal{K}(\mathcal{E},0)]^* \\
 = &\{M \in \mathbb{S}^n \mid \langle M, N \rangle \geq 0, \forall N \in \mathcal{K}(\mathcal{E},0)\} \\
 = &\left\{M \in \mathbb{S}^n \mid \left\langle M, \sum_{k=1}^p E_{\mathcal{C}_k}^\tr N_k E_{\mathcal{C}_k}\right\rangle \geq 0, \forall N_k \in K_k \right\} \\
  = &\left\{M \in \mathbb{S}^n \mid \sum_{k=1}^p\left\langle E_{\mathcal{C}_k}ME_{\mathcal{C}_k}^\tr,  N_k \right\rangle\geq 0, \forall N_k \in K_k \right\} \\
  = &\left\{M \in \mathbb{S}^n \mid  E_{\mathcal{C}_k}ME_{\mathcal{C}_k}^\tr \in K_k^*, k = 1, \ldots, p \right\} \\
  = & \mathcal{K}^*(\mathcal{E},?),
\end{aligned}
$$
where the second to last equality used the following fact: if any $E_{\mathcal{C}_t}ME_{\mathcal{C}_t}^\tr \notin K_t^*$ for some $t$, we can choose $N_t \in K_t$ such that 
$$
    \langle E_{\mathcal{C}_t}ME_{\mathcal{C}_t}^\tr, N_t \rangle < 0.
$$
Now, by choosing $N_k = 0 \in K_k, k \neq t$, we have 
$$
    \sum_{k=1}^p\left\langle E_{\mathcal{C}_k}ME_{\mathcal{C}_k}^\tr,  N_k \right\rangle = \langle E_{\mathcal{C}_t}ME_{\mathcal{C}_t}^\tr, N_t \rangle < 0,
$$
which contradict line 4. Thus, we must have $E_{\mathcal{C}_k}ME_{\mathcal{C}_k}^\tr \in K_k^*, k = 1, \ldots, p$.

% Similarly, we verify that 
% $$
% \begin{aligned}
%   &[\mathcal{K}(\mathcal{E},?)]^* \\
%  = &\{M \in \mathbb{S}^n \mid \langle M, N \rangle \geq 0, \forall N \in \mathcal{K}(\mathcal{E},?)\} \\
%   = &\{M \in \mathbb{S}^n \mid \langle M, E_{\mathcal{C}_k}^\tr N_k E_{\mathcal{C}_k} \rangle \geq 0, \forall N_k \in K_k\} \\
%   = &\{M \in \mathbb{S}^n \mid \langle E_{\mathcal{C}_k} M E_{\mathcal{C}_k}^\tr,  N_k  \rangle \geq 0, \forall N_k \in K_k\} \\
%     % = &\{M \in \mathbb{S}^n \mid \langle E_{\mathcal{C}_k} M E_{\mathcal{C}_k}^\tr,  N_k  \rangle \geq 0, \forall N_k \in K_k\} \\
%  ... 
% \end{aligned}
% $$
\end{proof}

For the containment in Equation \eqref{eq:dual_completion}:
$$
    [\mathcal{K}(\mathcal{E},?)]^* = \mathcal{K}^*(\mathcal{E},0).
$$

\begin{proof}
    Given any $M \in \mathcal{K}^*(\mathcal{E},0)$, we show $M \in [\mathcal{K}(\mathcal{E},?)]^*$. By definition (6), there exist $M_k \in K_k^*, k = 1, \ldots, p$, such that $$M = \sum_{k=1}^p E_{\mathcal{C}_k}^\tr M_kE_{\mathcal{C}_k}.$$
    We now verify that $\forall N \in \mathcal{K}(\mathcal{E},?)$
    $$
    \begin{aligned}
    \langle M, N \rangle &= \left\langle \sum_{k=1}^p E_{\mathcal{C}_k}^\tr M_k E_{\mathcal{C}_k}, N \right\rangle \\ & = \sum_{k=1}^p \left\langle   M_k, E_{\mathcal{C}_k}NE_{\mathcal{C}_k}^\tr  \right\rangle \\
    &= \sum_{k=1}^p \inp{M_k \in K^*_k}{N_k \in K_k} \quad \forall N_k \in K_k
    \\& \geq 0,\\
    \end{aligned}
    $$
    where the last inequality used the definition of $N \in \mathcal{K}(\mathcal{E},?)$.
    Thus, $M \in [\mathcal{K}(\mathcal{E},?)]^*$. We can conclude $\mathcal{K}^*(\mathcal{E},0) \subseteq [\mathcal{K}(\mathcal{E},?)]^*$
\end{proof}

The reverse inclusion ($\mathcal{K}^*(\mathcal{E},0) \supseteq [\mathcal{K}(\mathcal{E},?)]^*$) is only guaranteed to hold if the clique cover $\cs_k$ of $\es$ are disjoint: $\cs_k \cap \cs_{k'} = \empty \ \forall k \neq k'$.

\section{Applications to semidefinite optimization}

\label{sec:semidefinite} 

In this section, we develop inner and outer approximations of semidefinite programs using the notion of decomposed structured subsets, and discuss an application of $\mathcal{H}_{\infty}$ norm estimation of network systems. All code is publicly available at \url{https://github.com/zhengy09/SDPfw} within the folder \texttt{decomposed\_structured\_subsets}.

\subsection{Decomposed structured subsets in semidefinite programs}

A semidefinite program in primal form \eqref{Eq:SDPprimal} with $(X\in \psd_+^n)$ and dual form \eqref{Eq:SDPdual} $(Z \in \psd_+^n)$ will have matching optima $p^*= d^*$ when strong duality holds. By complementary slackness, $\inp{X}{Z} = 0$. Assume this semidefinite program has an aggregate sparsity pattern $\es$. With an optimization problem \eqref{Eq:SDPprimal} over $X \in \psd^n_+(\es, ?)$ and a cone set $\ks$, an upper bound is attained by imposing $X \in \ks(\es, ?)$ in \eqref{Eq:SDPprimal}, and a lower bound is found by restricting $Z \in \ks(\es, 0)$ in \eqref{Eq:SDPdual}. %In the same way, a problem posed over $\psd^n_+(\es, 0)$ will have an upper bound found by constraining $X \in \ks(\es, 0)$ and a lower bound by $Z \in \ks(\es, ?)$.

The conic optimization problem \eqref{Eq:SDPprimal} for a decomposed structured subset $\ks(\es, ?)$ is:
\begin{equation} \label{Eq:CliqueDependent}
    \begin{aligned}
    \min_{X} \quad & \langle C,X \rangle \\
    \text{subject to} \quad & \langle A_i,X \rangle = b_i, \; i=1, 2, \ldots, m \\
    & E_{\cs_k} X E_{\cs_k}^\tr  \in K_k.
    \end{aligned} 
\end{equation}

%Multiple forms of structure may be found in a single problem, such as when the SDP has symmetry and sparsity. The resultant optima will be closest to $p_{\text{SDP}}$ if all possible decompositions are applied before using structured subsets. Decompositions are not commutative; switching decomposition orders will generally result in different $\ks$'s. 

% As an example, consider a randomly generated block-arrow SDP with aggregate sparsity pattern $\es$ shown in blue in Figure \ref{fig:block_arrow}. This system has 80 equality constraints and a semidefinite block of size 160. The 15 blocks and the arrowhead are each of size 10, so each clique cone is $\psd^{20}_+$. The magenta entries form a coarser chordal completion $\es_F$ with larger clique cones (sizes 30, 40, 50). Table \ref{tab:block_arrow} shows costs of optimization over decomposed structured subsets. Cones $B_q$ are shorthand for block factor-width 2 matrices where each block has approximately $q$ elements (see \citep{zheng2019block} for further detail, $B_1 = \sdd$). By Theorem~\ref{T:AglerTheorem}, all entries of $K=\psd_+$ have the same optima. All entries $K = \dd$ are infeasible, and objectives decrease towards the bottom right corner of the table.

An example of decomposed structured subsets in action is a random SDP with block arrow sparsity pattern with 80 equality constraints, where each of the 15 blocks has size 10 and the arrowhead has width 10 (see Figure \ref{fig:block_arrow} for the aggregate sparsity pattern). The original SDP has $X \in \psd_+^{160}$, and chordal decomposition has $X_k \in \psd_+^{20}$ that are equal in the $10 \times 10$ bottom right corner (blue pattern $\es$). A coarser chordal decomposition is the union of blue and magenta blocks in Figure \ref{fig:block_arrow} (fill-in $\es_F$), which has block sizes $(\psd_+^{50})^2  \times \psd_+^{41} \times (\psd_+^{30})^2$ that are each equal in the bottom right corner. Clique consistency for $\es_F$ adds 230 equality constraints, and for $\es$ adds 770. 

Cost values of this SDP and its approximants are recorded in Table \ref{tab:block_arrow}. Rows are different structured subsets and the columns apply the structured subsets: imposing that $X \in K$ and then that the cliques of $X_k \in K_k$, where cliques are set based on the graphs $\es_F$ and $\es$. We introduce shorthand $B_k$ as a cone of block factor-width 2 matrices where each block has $k$ components (so $B_1$ = $\sdd$) and membership constraints in $\psd^{2k}_+$ are imposed. As an example, the row-column pair $(B_5, K)$ corresponds to the cone $B_5^{160}$ and $(B_10, K(\es_F, ?)$ refers to the cone $B_10^{160}(\es_F, ?)$.'

By Agler's theorem, all entries of $K=\psd_+^{\new{160}}$ have the same optima. All entries $K = \dd^{\new{160}}$ are infeasible, and objectives decrease towards the bottom right corner of the table as expected in the above containment analysis.
Table \ref{tab:block_arrow} demonstrates that merging blocks together may degrade the resultant approximation quality. Merging strategies such as \cite{garstka2020clique} must therefore be used with caution. Even well-chosen merges that speed up program execution may worsen the approximated SDP bound.

\begin{table}[ht]
  \begin{minipage}[b]{0.35\linewidth}
    \centering
    \includegraphics[width=0.6\linewidth]{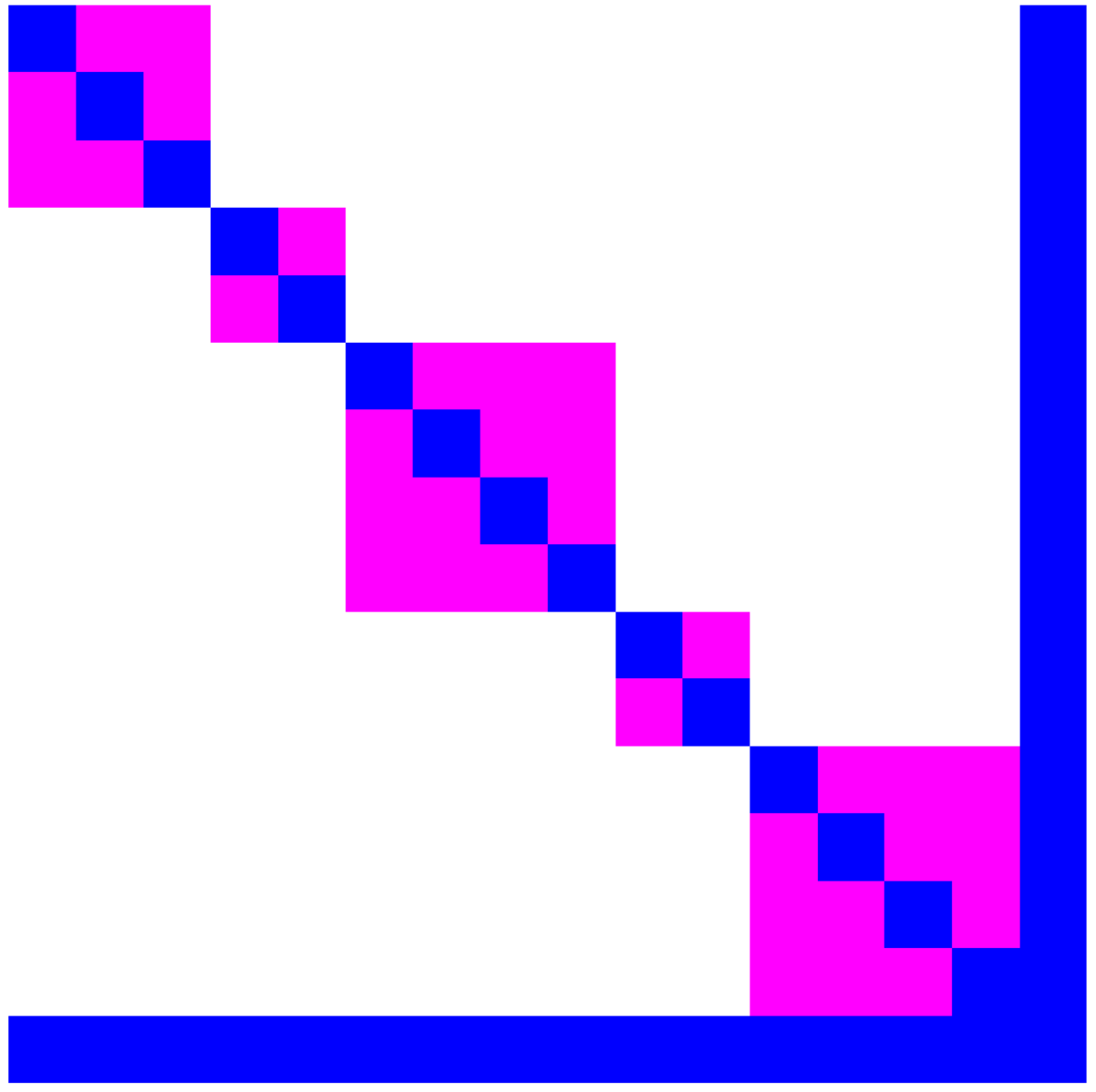}
    \captionof{figure}{\label{fig:block_arrow} Block arrow sparsity pattern and extension}
    \label{fig:image}
  \end{minipage}
 \hfill
  \begin{minipage}[b]{0.6\linewidth}
  \setlength{\abovecaptionskip}{0mm}
  \setlength{\belowcaptionskip}{0mm}
    \renewcommand\arraystretch{1.0}
    \centering
    \begin{tabular}{ llll }
 \toprule 
  & $K$ & $K(\es_F,?)$ & $K(\es, ?)$ \\
  \hline
$\dd$ & Inf. & Inf. &Inf.\\ 
$B_1$ & 64.5 & 34.7 & 19.4\\ 
$B_2$  &51.4 & 27.1 & 13.9\\ 
$B_5$  &32.1 & 15.0 & 5.34\\ 
$B_{10}$  &20.8 & 7.10 & -1.23\\ 
$\psd_+$  &-1.23 & -1.23 & -1.23 \\
\bottomrule
    \end{tabular}
    \caption{\label{tab:block_arrow} Cost vs. subset and pattern in Fig. \ref{fig:block_arrow}}
  \end{minipage}%
\end{table}

% \begin{figure}[ht]
% \centering
% % \begin{minipage}{0.39\linewidth}

% % \centering
% \includegraphics[width=0.4\linewidth]{img/block_arrow_sparsity_15.png}
% \caption{\label{fig:block_arrow} Block arrow sparsity pattern and extension}
% \end{figure}
% % \caption{Block Arrow Pattern}\label{fig:block_arrow}
% % \end{minipage}
% % \begin{minipage}{0.59\linewidth}
% \begin{table}
% \setlength{\abovecaptionskip}{0mm}
%   \setlength{\belowcaptionskip}{0mm}
%     \renewcommand\arraystretch{1.0}
% \caption{\label{tab:block_arrow} Cost vs. subset and pattern in Figure \ref{fig:block_arrow}}
% \centering
%     \begin{tabular}{llll}
%     \toprule 
%   & $K$ & $K(\es_F,?)$ & $K(\es, ?)$ \\
%   \hline
% $\dd$ & Inf. & Inf. &Inf.\\ 
% $B_1$ & 64.5 & 34.7 & 19.4\\ 
% $B_2$  &51.4 & 27.1 & 13.9\\ 
% $B_5$  &32.1 & 15.0 & 5.34\\ 
% $B_{10}$  &20.8 & 7.10 & -1.23\\ 
% $\psd_+$  &-1.23 & -1.23 & -1.23 \\
% \bottomrule
% \end{tabular}
% % \end{minipage}
% \end{table}

\subsection{Certifying Optimality}

\label{sec:certify_optimum}

Karush-Kuhn-Tucker (KKT) conditions can be used to certify if an SDP approximated by structured subsets reaches the same cost as the original SDP \cite{ahmadi2017optimization}. For a cone $K \subset \psd_+ \subset K^*$ and an SDP of type \eqref{Eq:SDPprimal}, KKT relations on the left side will hold for an optimal primal-dual triple $(X, y, Z) \in (K^n, \R^m, K^{n^*})$ while optimizing over $X \in K$. All KKT conditions over $K$ in Equation \eqref{eq:kkt_k} are also guaranteed to hold over $\psd_+$ in Equation \eqref{eq:kkt_psd} except for dual feasibility ($Z \in K^*$):

\noindent\begin{minipage}{.5\linewidth}
  \begin{align}
  \label{eq:kkt_k}
    X &\in K \nonumber \\
    Z &\in K^* \\
    XZ &= 0 \nonumber \\
    b_i &= \inp{A_i}{X} \mid_{i=1}^m \nonumber \\
    C &= S + \sum_{i=1}^m y_i A_i \nonumber,
  \end{align}
\end{minipage}%
\begin{minipage}{.5\linewidth}
  \begin{align}
  \label{eq:kkt_psd}
    X &\in \psd_+ \nonumber \\
    Z &\in \psd_+ \\
    XZ &= 0 \nonumber \\
    b_i &= \inp{A_i}{X} \mid_{i=1}^m \nonumber \\
    C &= S + \sum_{i=1}^m y_i A_i \nonumber.
  \end{align}
\end{minipage}

% \begin{multicols}{2}
% \noindent 
%   \begin{align}
%     X &\in K \nonumber \\
%     Z &\in K^* \\
%     XZ &= 0 \nonumber \\
%     b_i &= \inp{A_i}{X} \mid_{i=1}^m \nonumber \\
%     C &= S + \sum_{i=1}^m y_i A_i \nonumber 
%   \end{align}
%   \begin{align}
%     X &\in \psd_+ \nonumber \\
%     Z &\in \psd_+ \\
%     XZ &= 0 \nonumber \\
%     b_i &= \inp{A_i}{X} \mid_{i=1}^m \nonumber \\
%     C &= S + \sum_{i=1}^m y_i A_i \nonumber 
%   \end{align}
% \end{multicols}

If the dual matrix $Z$ is PSD, then $(X, y, Z)$ solves the original SDP with the same optimal cost.
Checking if $Z$ has a negative eigenvalue can be accomplished by inverse iteration (power method on $Z^{-1}$).
SDP-optimality of decomposed structured subsets can be certified in the same framework given a set of clique cones $\ks$. When finding lower bounds to semidefinite programs, the clique cones $\ks$ have $K_k \supseteq \psd_+$. Tightness is certified if each $X_k \in \psd_+$. Upper bounds have $K_k \subseteq \psd_+$. For each clique $\cs_k$ in the clique cone $\ks$, check if the corresponding dual block $Z_k \in \psd_+$. The dual clique blocks $Z_k$ can be obtained by computing $Z = C - \sum_{i=1}^m {y_i A_i}$.

% For a problem in $\ks(\es, ?)$, all cliques are in $\ks$ and they agree on overlaps to complete into a PSD matrix. 

% The dual $Z \in \ks(\es, ?)^* = \ks^*(\es, 0)$ will be a sparse matrix, and inverse iteration can occur efficiently through sparse linear system algorithms. On the sparse cone $\ks(\es, 0)$, the dual solution $Z \in \ks(\es, ?)$ will agree on overlaps, and all cliques in $Z$ must be checked to ensure that $Z_k\in \psd_+$. In the practical case where decomposition methods introduce e

% In the case of $K = \sdd$, $K^* = \sdd^*$ is the set of symmetric matrices with all $2\times2$ principal minors PSD. If $Z \in \sdd^*$ is additionall

% \caption{\label{tab:block_arrow} Cost vs. (decomposed) structured subset}
%\end{minipage}

% \begin{figure}[ht]
% \centering
% \begin{minipage}{0.39\linewidth}
% \centering
% \includegraphics[width=\linewidth]{img/block_arrow_sparsity_15.png}
% \caption{Block Arrow Sparsity Pattern}\label{fig:block_arrow}
% \end{minipage}
% \begin{minipage}{0.59\linewidth}
% \centering
% \begin{tabular}{l|lll}
%   & $K$ & $K(\es_F,?)$ & $K(\es, ?)$ \\
% $\dd$ & Inf. & Inf. &Inf.\\ 
% $B_1$ & 64.5 & 34.7 & 19.4\\ 
% $B_2$  &51.4 & 27.1 & 13.9\\ 
% $B_5$  &32.1 & 15.0 & 5.34\\ 
% $B_{10}$  &20.8 & 7.10 & -1.23\\ 
% $\psd_+$  &-1.23 & -1.23 & -1.23
% \end{tabular}
% \caption{\label{tab:block_arrow} Cost vs. (decomposed) structured subset}
% % \caption{\label{tab:block_arrow} Cost vs. (decomposed) structured subset}
% \end{minipage}
% \end{figure}

\subsection{Decomposed Change of Basis}

Decomposed structured subsets are compatible with the change of basis algorithm as reviewed in Section \ref{Sec:ChangeofBasis}. 
% Assume that Problem \eqref{Eq:CliqueDependent} is solved to find a solution $X_0 \in \ks(\es, ?)$. 
\new{Assume that $X_0 \in \ks(\es, ?)$ is a solution to Problem \eqref{Eq:SDPprimal}}. Define Cholesky factorization matrices $L_k^0$ for each clique $k = 1 \ldots p$ such that \new{$L_k^0 L_k^{0\tr} = E_{\cs_k} X_0 E_{\cs_k}^\tr$.}
% \begin{equation}
% \label{eq:change_factor}
%     L_k^0 L_k^{0\tr} = E_{\cs_k} X_0 E_{\cs_k}^\tr\new{.}
% \end{equation}
% It is not necessary for $X_0$ to have been PSD-completed to find Cholesky factorization matrices $L_k$, the matrices $L_k$ are purely a property of entries supported on clique $\cs_k$.
The next iteration of the change of basis algorithm will solve
\begin{equation} \label{Eq:CliqueDependent_changed}
    \begin{aligned}
   X_1 = \arg\min_{X} \quad & \langle C,X \rangle \\
    \text{s.t.} \quad & \langle A_i,X \rangle = b_i, \; i=1, 2, \ldots, m \\
    & E_{\cs_k} X E_{\cs_k}^\tr   \in K_k(L_k^0). 
    \end{aligned} 
\end{equation}
% The cone constraint in \eqref{Eq:CliqueDependent_changed} for $X$ may also be expressed as $X \in \ks(\mathcal{L}^0)(\es, ?)$ for a set of basis change matrices $\mathcal{L}^0 = \{L_k\}_{k=1}^p$. 
The solution to \eqref{Eq:CliqueDependent_changed} can be used to find a new set of factor matrices $L^1_k$ \new{ by finding $L_k^1 L_k^{1\tr} = E_{\cs_k} X_1 E_{\cs_k}^\tr$.}
% by the method in \eqref{eq:change_factor}. 
% In this manner, an accumulated basis change is formed between iterations as $L_k = L^t_k L^{t-1}_k \ldots L^1_k L^0_k$ up to iteration $t$. 
Each clique $\cs_k$ is described by basis $L_k$, and different bases may describe the same elements of $X$ on clique-overlaps. 

\begin{remark}
    \new{Performing a decomposed change-of-basis over $K^n(\es, ?)$ will result in a lower cost as compared to applying change-of-basis over $K^n$ at the first iteration. No conclusions can be drawn after the first iterations.
% \end{remark}
% \begin{proof}
% At the first iteration, $K^n \subseteq K^n(\es, ?)$, and therefore the cost from $K^n$ will have a value $\geq$ the cost from $K^n(\es, ?)$. Basis change matrices $L$ for $\{X \in K^n(L)\}$ and $\{L_k\}_{k=1}^p$ for $\{X_k \in K_k^{\abs{\cs_k}}(L_k)\}$ are computed to start the second change-of-basis iteration. 
% These basis-changed cones are no longer related, and no inferences can be drawn about long-term cost behavior. 
In experiments, the cost sequence obtained from performing change-of-basis over $K^n$ remains above $K^n(\es, ?)$'s cost sequence. }
% This behavior is observed in Figure \ref{fig:change_of_basis_arrow}.
% \end{proof}
\end{remark}

\begin{figure}[ht]   
\centering   
\begin{subfigure}[t]{0.35\linewidth}        \includegraphics[width=\linewidth]{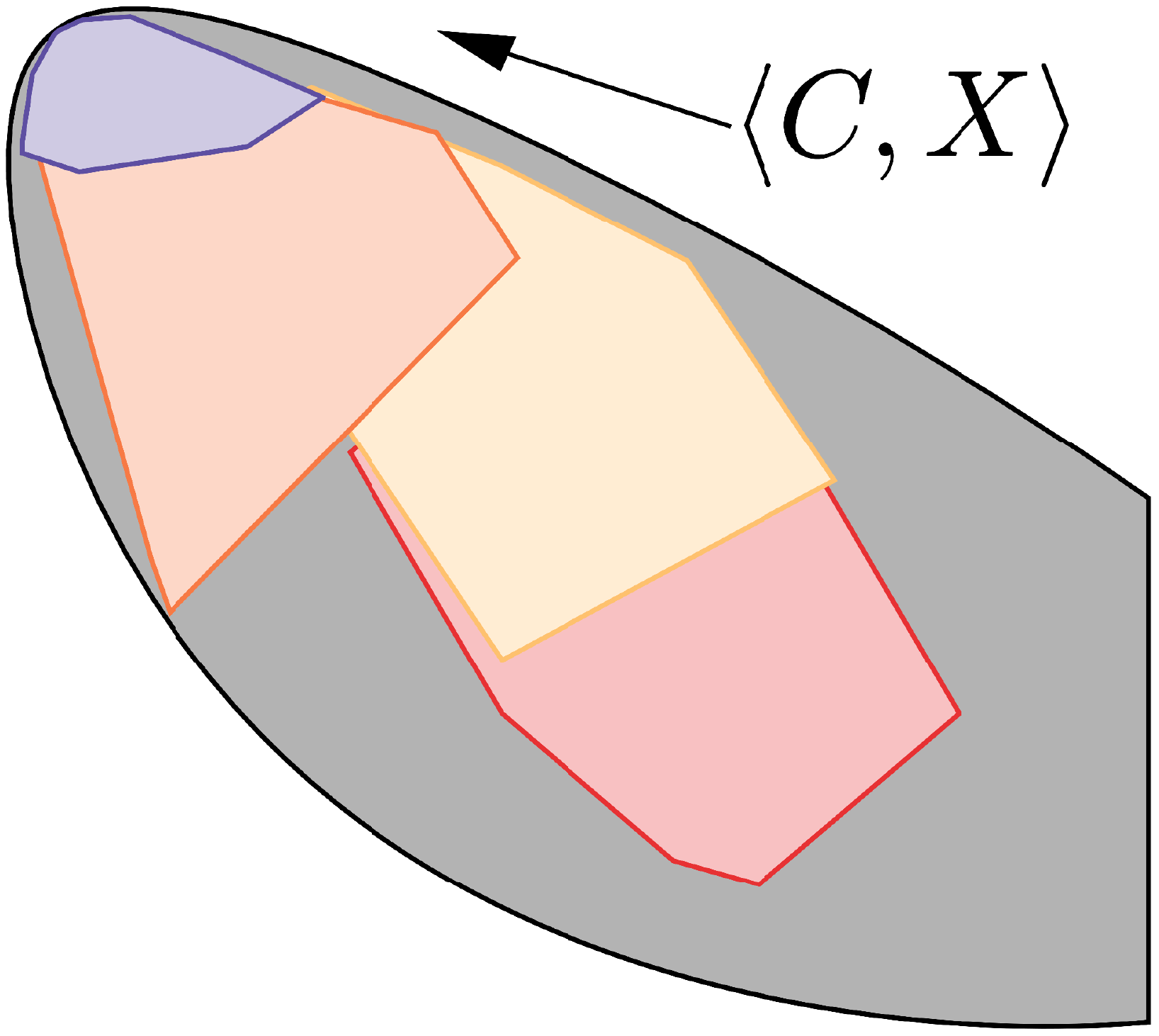}       
%\caption*{Start from ${\color{dd} \dd}$}
\label{fig:change_of_basis_std}    
\end{subfigure}   
\qquad 
\begin{subfigure}[t]{0.35\linewidth}        \includegraphics[width=\linewidth]{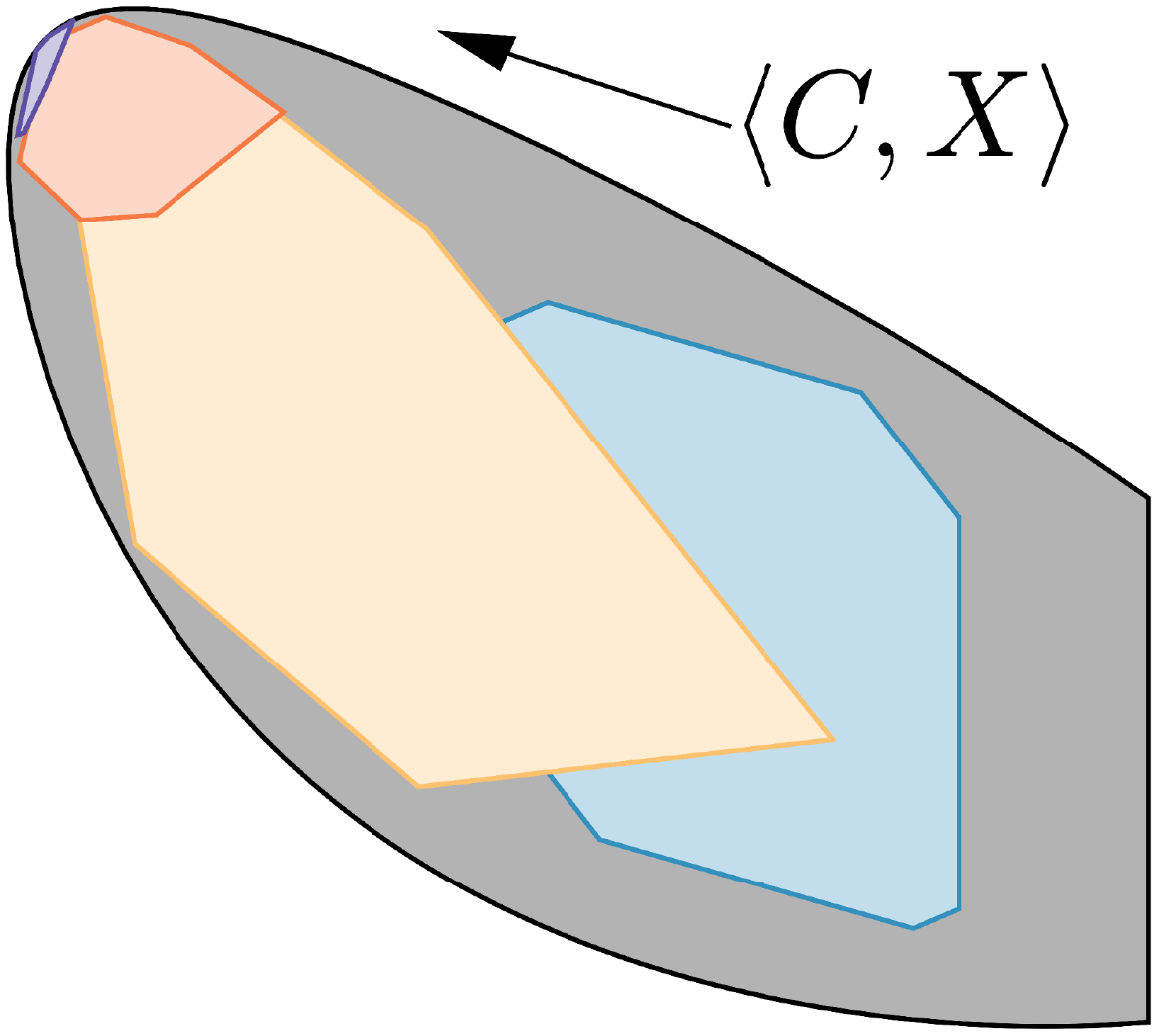}        
%\caption*{Start from ${\color{cdd} \dd(\es, ?)}$}     
\label{fig:change_of_basis_split}    \end{subfigure}    
\caption{\label{fig:change_of_basis} Decomposed vs. Standard Change of Basis on \eqref{eq:exampleM}. Left: Start from ${\color{dd} \dd^4}$; Right: Start from ${\color{cdd} \dd^4(\es, ?)}$. }
\end{figure}

Figure \ref{fig:change_of_basis} illustrates the change of basis technique on a standard (left) and decomposed (right) structured subsets in the direction of $\inp{C}{X}$. The intermediate costs are recorded in Table \ref{tab:change_of_basis} for the first three iterations:

\begin{table}[ht]
\setlength{\abovecaptionskip}{0mm}
   \setlength{\belowcaptionskip}{0mm}
    \renewcommand\arraystretch{1.0}
\centering
\caption{\label{tab:change_of_basis} Decreasing costs over Change of Basis}
\begin{tabular}{lcccc}
\toprule
              & \multicolumn{4}{c}{Change of Basis Iteration} \\
              \hline
              & 0     & 1     & 2     & 3     \\
$\dd^4$         & -1.41 & -2.50 & -3.08 & -3.15 \\
$\dd^4(\es, ?)$ & -1.41 & -3.02 & -3.13 & -3.17 \\
\bottomrule
\end{tabular}
\end{table}

\begin{figure}[ht]
    \centering
    \includegraphics[width=0.6\linewidth]{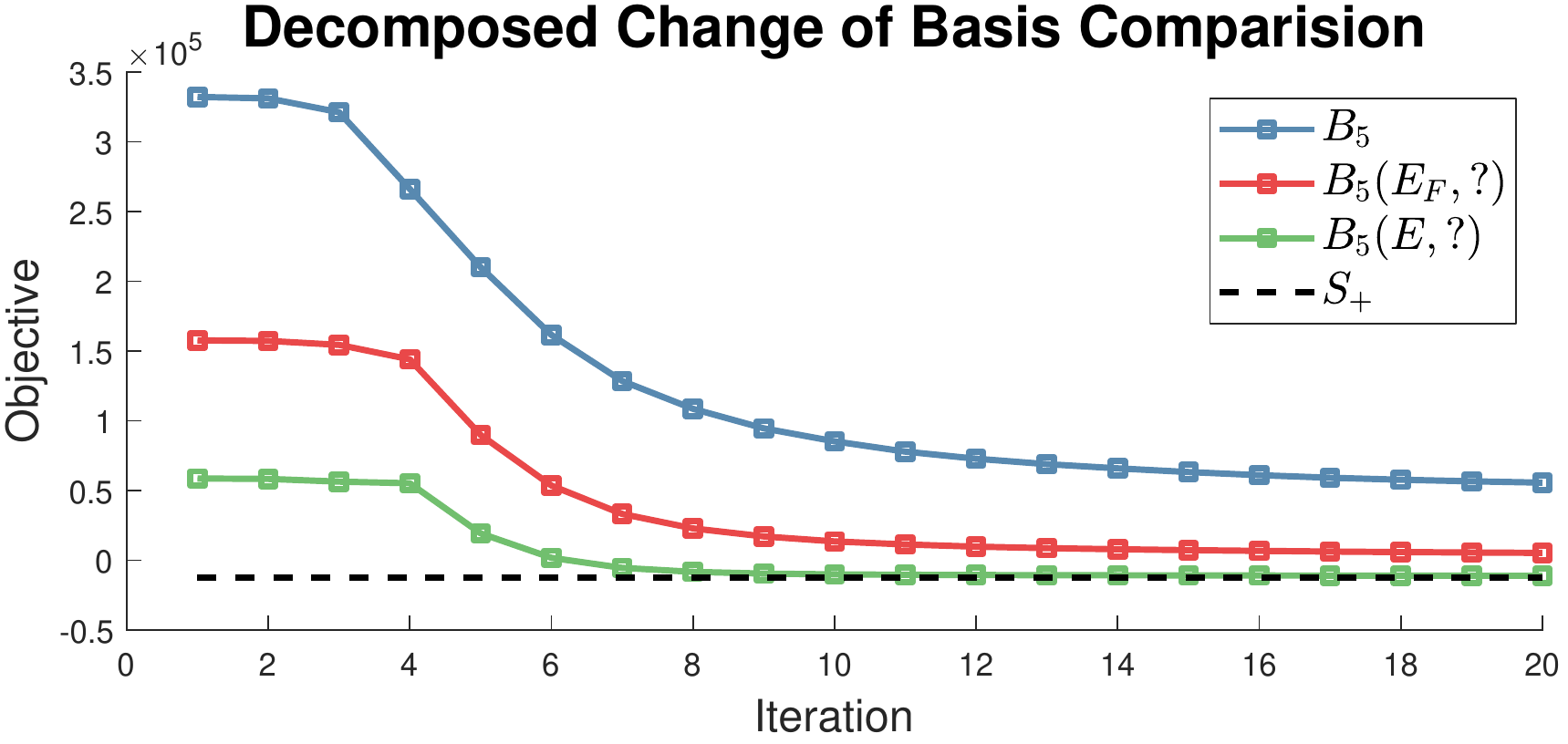}
    \caption{Change of basis on Block Arrow SDP}
    \label{fig:change_of_basis_arrow}
\end{figure}

Figure \ref{fig:change_of_basis_arrow} shows the output of the change of basis algorithm for the cone $B_5^{\new{160}}$ on the block arrow system shown in Figure \ref{fig:block_arrow}. Over the course of 20 iterations, the basis-changed cone starting with $B_5^{\new{160}}(\es, ?)$ (green curve) \new{eventually} matches the SDP optimum.

A similar process can be done over the sparse cone $Z \in \ks(\es, 0)$ (dual SDP), where bases $\mathcal{L} = \{L_k\}_{k=1}^p$ are tracked for each clique component $Z_k$ forming the clique-sum $Z = \sum_k{E_{\cs_k}^\tr Z_k E_{\cs_k}}$ for $Z \in \ks(\mathcal{L})(\es, 0)$.

\subsection{H-infinity Norm Estimation for Networked Systems}
\label{sec:hinf}
Here, we present a special applications of SDPs in $\mathcal{H}_\infty$ norm estimation. Consider a state-space stable dynamical system $G(s)$:
\begin{align*}
    \dot{x} &= A x + B u, \\
    y &= C x + D u.
\end{align*}
The $\mathcal{H}_{\infty}$ norm of $G(s)$ is the supremum over frequencies $\omega$ of the maximum singular value of $G(j \omega)$.
The norm $\norm{G}_{\infty}$ is finite when $A$ is Hurwitz. 
The Bounded Real Lemma can be used to find upper bounds on $\norm{G}_\infty$:
 \begin{theorem} [Bounded Real Lemma \citep{boyd1994linear}] \label{T:bounded_real}
    The following statements are equivalent:
    \begin{enumerate}
        \item $\norm{G}_{\infty} < \gamma,$
        \item There exists a $P \succ 0$ such that 
        \[\begin{bmatrix}P A + A^\tr P + C^\tr C& P^\tr B + C^\tr D\\ B^\tr P + D^\tr C & -\gamma^2 I \end{bmatrix} \prec  0. \]
    \end{enumerate}
     \end{theorem}

% \begin{align*}
% \norm{G(s)}_{\infty} &= \min_{\gamma, P} \gamma \\
% &P \geq 0
% &\begin{bmatrix}P A + A^\tr P & P^\tr B & C^\tr \\ B^\tr P & -\gamma I & D^\tr \\ C & D & -\gamma I \end{bmatrix} \geq 0 
% \end{align*}

%In the case of a networked system, such a $P$ matrix 

% {\color{red} Introduce networked systems, restrictions to block diagonal $P$, and how that does/does not add conservatism to system}.

%$\norm{G}_{\infty}$ can be estimated by minimizing $\gamma^2$ subject to $P \succ 0$ and the Bounded Real Lemma holding. 
If the dynamical system is sparse (has a network structure), a dense $P \succ 0$ 
%Computational efficiency can be leveraged by using chordal sparsity if the system has a network structure. 
%Using a dense matrix $P \geq 0$ 
will give the tightest $H_\infty$ approximation but will destroy the sparsity pattern. Choosing a $P$ structure to be compatible with the LMI sparsity pattern will form a computationally tractable upper bound of $\norm{G(s)}_\infty$. 
%A block-diagonal $P$ such that the size of each agent's block is equal to its number of states respects the network sparsity pattern 
One structure on a block-diagonal $P$ that respects the network sparsity pattern is when the size of each agent's block in $P$ equals its number of states 
\citep{zheng2018scalable}. %\cite{Sootla2019OnTE} has explored cases under which a block-diagonal Lyapunov function exists, which indicates a tight $H_\infty$ estimation.
%Network systems: help? \cite{zheng2018scalable, Sootla2019OnTE}
%\cite{Sootla2019OnTE} present ci

\begin{figure}[hb]
    \centering
    \includegraphics[width=0.7\linewidth]{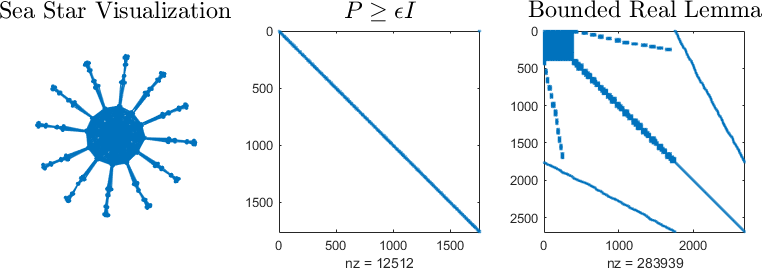}
    \caption{Sea Star network topology and LMI sparsity}
    \label{fig:sea_star_network}
\end{figure}

As an example of applying decomposed structured subsets to $\mathcal{H}_{\infty}$ estimation, we present the `sea star' networked system. The sea star system is composed of a set of agents clustered into a head and a set of arms. Each agent has internal linear dynamics ($n_i$ states, $m_i$ inputs, $d_i$ outputs), and they communicate and respond to a sparse selection of other agents. The left panel of Figure \ref{fig:sea_star_network} shows a sea star network with 70 densely connected agents in the head and other agents distributed into 12 arms. Each arm is composed of 2 densely connected `knuckles'. Each knuckle has 10 agents, and every knuckle in the arm communicates with 4 agents in the next and previous knuckle (or the head as appropriate). The individual agent dynamics combine to form global dynamics $[A, B, C, D]$, where $A$ is Hurwitz. 

% {\color{red} Change the parameters as appropriate. Mosek runs out of memory for this system. Lowered head from 60 to 41 then to 30, may keep dropping due to memory constraints}.

% \begin{figure}[ht]
%     \centering
%     \includegraphics[width=0.8\linewidth]{img/network_plot-1.png}
%     \caption{Sea Star network topology and adjacency matrix }
%     \label{fig:sea_star_network}
% \end{figure}

Estimating $\norm{G(s)}_{\infty} = \norm{C(sI-A)^{-1}B + D}_{\infty}$ can be accomplished by using the bounded real lemma to minimize $\gamma^2$. The resultant LMI has two semidefinite variables, and the center and right panels of Figure \ref{fig:sea_star_network} displays the sparsity pattern of constraints on these variables. The top left corner of the Bounded Real LMI shows a structure induced by the network interconnections. On their own, the two semidefinite blocks are of size 1760 and 2691.

% \begin{figure}[ht]
%     \centering
%     \includegraphics[width=0.8\linewidth]{img/lmi-1.png}
%     \caption{Sparsity pattern of the $H_\infty$ LMI}
%     \label{fig:sea_star_sparsity}
% \end{figure}

\begin{figure}[ht]
    \centering
    \includegraphics[width=0.6\linewidth]{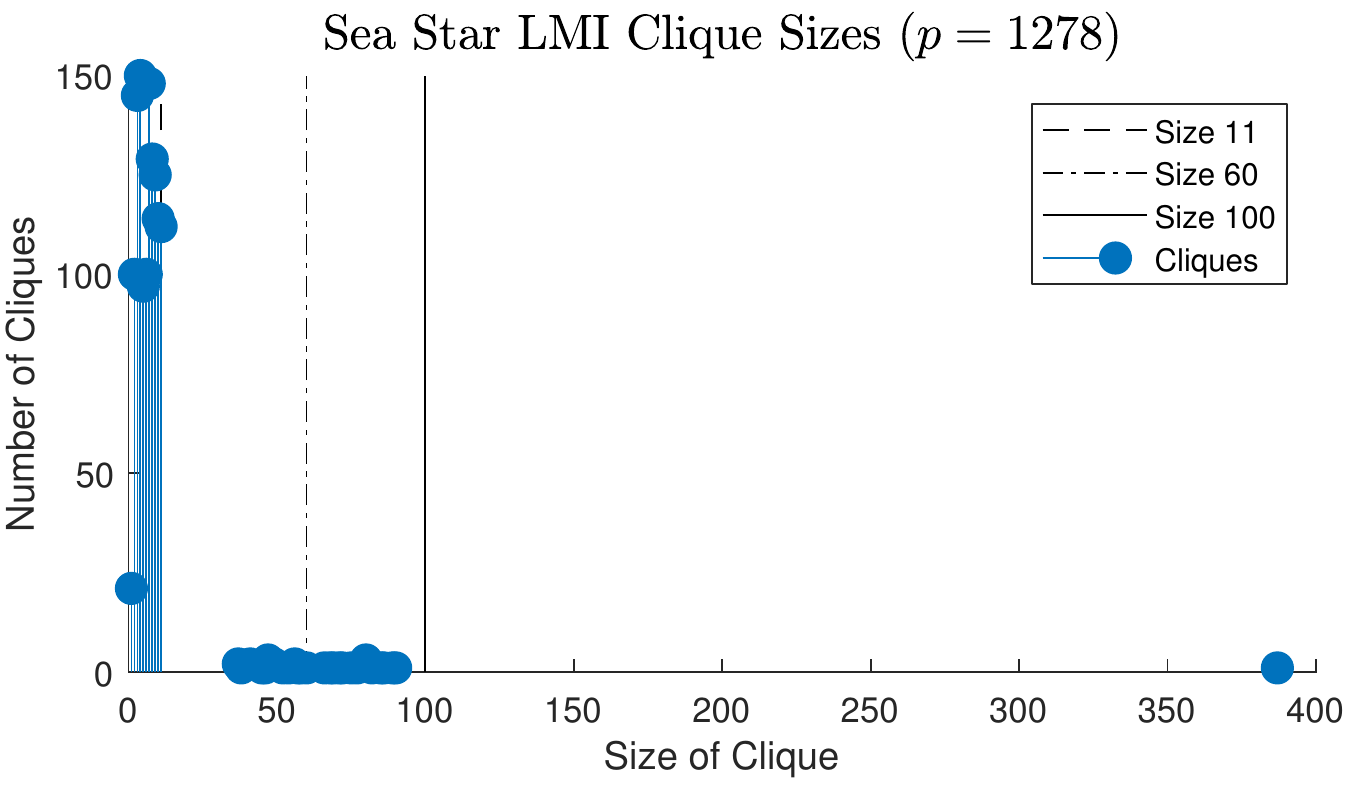}
    \caption{Maximal cliques of $\es$ from LMI with size thresholds}
    \label{fig:sea_star_clique}
\end{figure}

This LMI system strongly exhibits chordal sparsity with edges $\es$, and can be posed as an optimization problem over the cone $\psd_+(\es, 0)$. $\abs{\cs_k}$ is shown in Figure \ref{fig:sea_star_clique}. There is a run of cliques of sizes ranging from 1-11, a set of cliques from sizes 37-90, and a solitary clique of size 387. 

Results of $\mathcal{H}_{\infty}$ norm estimation of the sea star system are presented in Figures \ref{fig:sea_star_upper_time},  \ref{fig:sea_star_lower_gamma}, and \ref{fig:sea_star_lower_time}. Columns are cones $K$ where $K= \dd$ or $K=B_q$ if $q$ is an integer \ref{sec:factor_width}. Rows are size thresholds: the cone $K(\es, 0)$ has all cliques in $K$, and $\ks_{60}(\es, 0)$ is a mixed cone where cliques with $\abs{\cs} \leq 60$ are PSD and $\abs{\cs} > 60$ are in $K$. All experiments were written in MatlabR2018a and performed on Mosek \citep{andersen2000mosek} on a Intel i7 CPU with a clock frequency of 2.7GHz and 16.0 GB of RAM. %For reference, optimizing over $\psd_+$ (without taking into account sparsity) returns a cost of 1.224 in 3308 seconds.

%Run the experiments again, plot make plots in python 

\begin{table}[ht]
    \centering
    \caption{Time to find $\gamma$ by upper bound $K$ (minutes)}
    \label{fig:sea_star_upper_time}
        \includegraphics[width=0.8\linewidth]{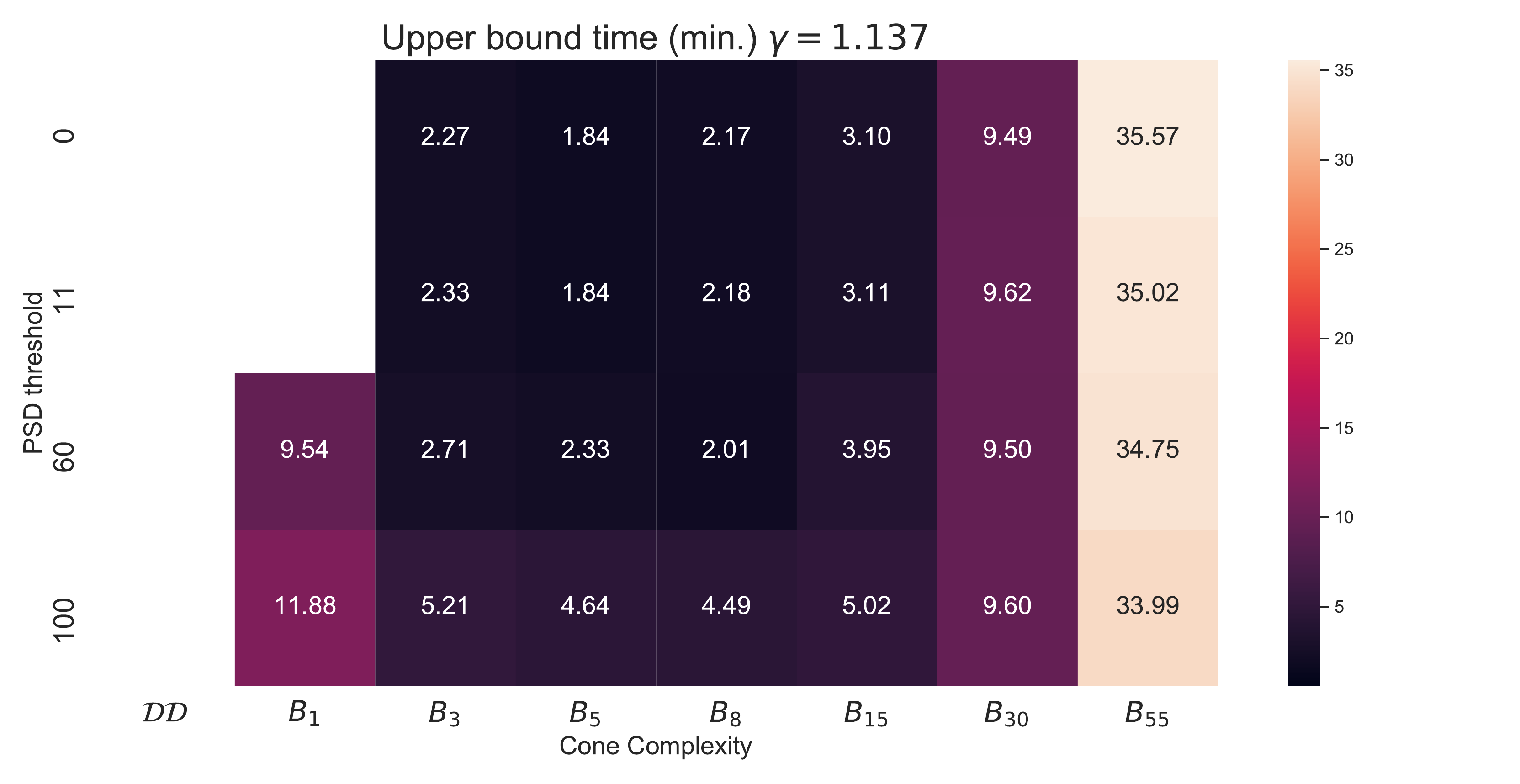}
\end{table}

Times in Figure \ref{fig:sea_star_upper_time} were measured solving the primal program over $K$. All displayed values achieved the SDP optimal solution, as certified in Section \ref{sec:certify_optimum}. 
The cones $\dd$ with size thresholds 0 and 11 were primal infeasible, other non-displayed values did not attain the optimal $\gamma = 1.137$. 
The cone $B_5$ was fastest at 1.84 minutes. 

\begin{table}[ht]
    \centering
    \caption{$\gamma$ found by lower bound $K^*$}
    \label{fig:sea_star_lower_gamma}
        \includegraphics[width=0.8\linewidth]{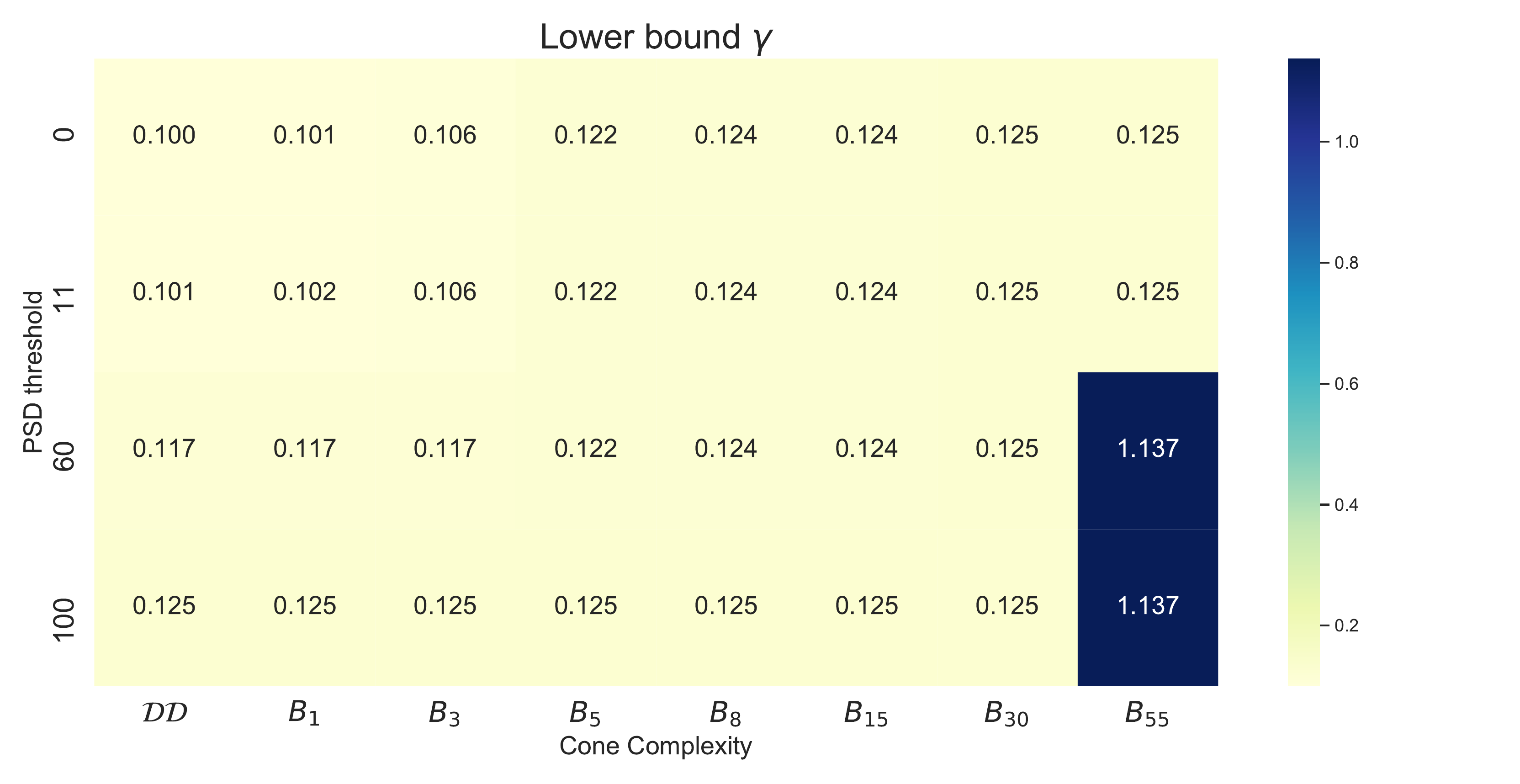}
\end{table}

\begin{table}[ht]
    \centering
    \caption{Time to find lower bound $\gamma$ over $K^*$ (minutes)}
    \label{fig:sea_star_lower_time}
        \includegraphics[width=\linewidth]{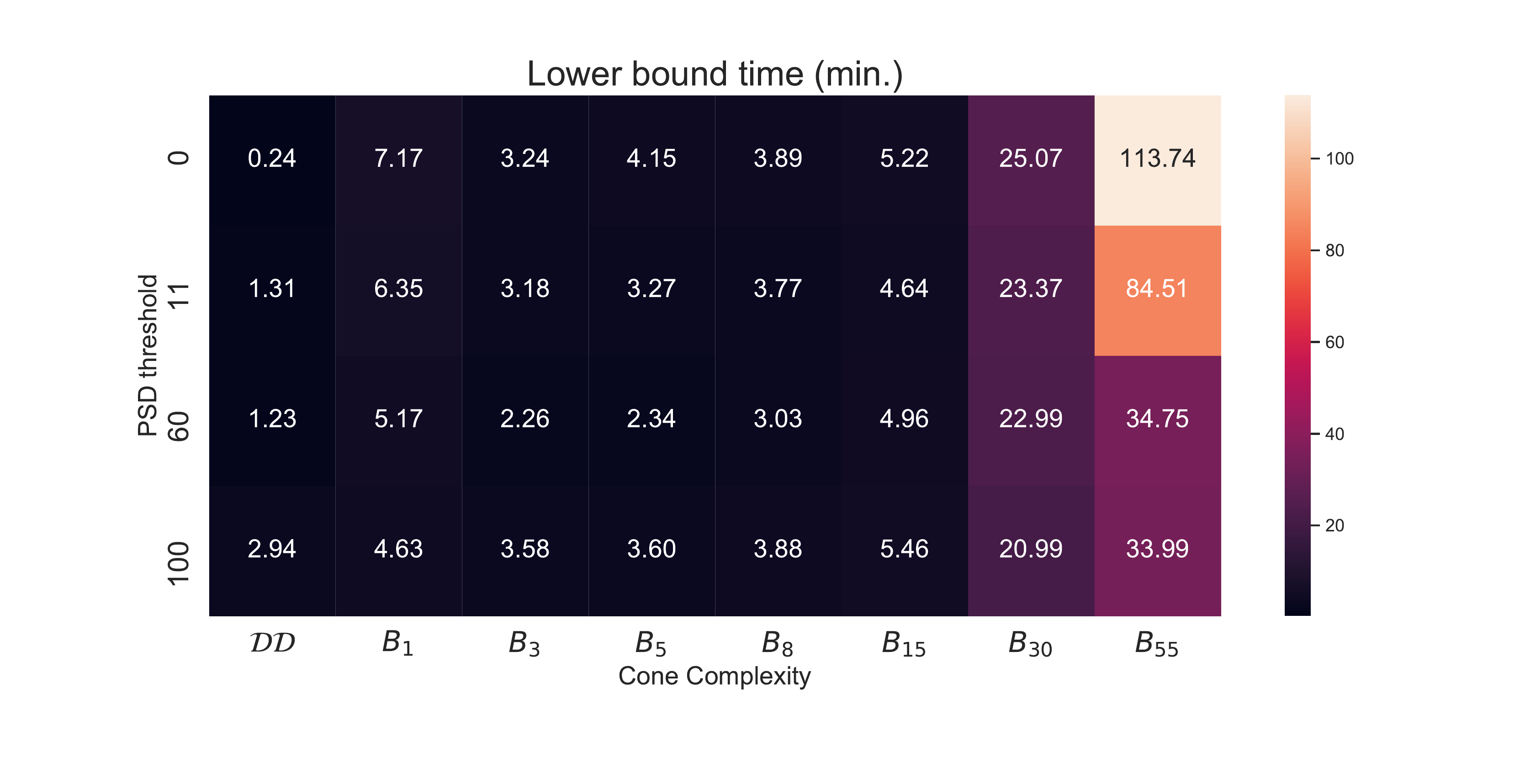}
\end{table}

Figure \ref{fig:sea_star_lower_gamma} displays lower bounds for $\gamma$ by over the dual cone $\ks^*(\es, 0)$. Lower bounds tighten as cone complexity and size thresholds increase. The true $\gamma$ is obtained with a block-size of 55 and size-thresholds of 60 and 100 taking 34.8 and 34.0 minutes. Figure \ref{fig:sea_star_lower_time} contains the time taken to find all lower bounds. The computer running experiments ran out of memory attempting to solve the LMI over $\psd_+(\es,0)$.

\section{Applications to Polynomial Optimization}

\label{sec:polynomial}
This section reviews sum-of-squares methods for approximating the optimal values of polynomial optimization problems, and demonstrates how decomposed structured subsets may be applied.

\subsection{Preliminaries for Polynomial Optimization}
A polynomial optimization problem may be approximated by semidefinite programming. Other methods include using nonsymmetrtic cone optimization \cite{papp2019sum}, applying Sketchy-CGAL \cite{yurtsever2021scalable} under limited memory requirements and low rank SDP solution structure, and exploiting the Constant Trace Property of moment relaxations \cite{mai2020hierarchy, mai2020exploiting}.

The task of minimizing a polynomial  $p(x) \in \R[x]_{\leq d}$ of bounded degree $d$ is equivalent to solving \citep{lasserre2010moments}:
\begin{equation}
    \label{eq:poly_nonneg}
    \begin{aligned}
        p^* =& \max_{\gamma} {\gamma} \\
        & p(x) - \gamma \geq 0.
    \end{aligned}
\end{equation}

%Let $\K$ be a compact semialgebraic region defined by inequality constraints $g_i(x) \leq0$ and equality constraints $h_j(x)=0$. Polynomial optimization problems to minimize $p(x)$ over $\K$ is equivalent to solving \citep{lasserre2010moments}:

% \begin{equation}
%     \begin{aligned}
%         p^* &= \min_{x}{p(x)} \\
%             & g_i(x) \leq0 \\
%             & h_j(x) =0
%     \end{aligned}
% \end{equation}

% If $\K$ is Archimedean (or compact with the redundant constraint $\norm{x}^2-R \leq0$ where the ball $B_R \supset \K$), then the Putinar Positivstellensatz forms an equivalent problem:

% \begin{equation}
%     \begin{aligned}
%         p^* &= \max_{\gamma} {\gamma} \\
%         &P(x) = p(x) - \gamma + \sum_i {\sigma_i(x)g_i(x)} + \sum_j{h_j(x) \phi_j(x)} \geq0 \\
%         &\sigma_i(x) \geq0 \qquad \phi_j(x) \in \R[x]
%     \end{aligned}
% \end{equation}

%Detecting or imposing polynomial non-negativity is an NP-hard nonconvex problem. Polynomial optimization problems may involve additional parameters: find $\theta \in \Theta$ such that $p(x; \theta) \geq0 \  \forall x$ 
The imposition of a polynomial nonnegativity constraint is generically NP hard. Sum-of-squares (SOS) methods offer a convex relaxation of polynomial nonnegativity through SDPs. The square of a real number is always nonnegative, so a polynomial $q(x) = \sum_{i} \sigma_i(x)^2$ for polynomoials $q(x) \R[x]$ is also nonnegative. The cone of SOS polynomials $q(x) \in \Sigma[x]_{\leq d}$ is the set of polynomials of degree $\leq d$ that admit such a decomposition into terms $q(x)$. If $v_d(x)$ denotes a monomial map lifting $x$ into a set monomials up to degree $d$, then every SOS polynomial $\sigma(x) \in \Sigma[x]_{\leq d}$ has an equivalent description of $\sigma(x) = v(x)^\tr Q v(x)$ for some 
$Q \succeq 0$. The SOS relaxation of degree $d$ of Equation \eqref{eq:poly_nonneg} is:
\begin{equation}
    \begin{aligned}
        \label{eq:poly_sos}
        p_d^* =& \max_{\gamma} {\gamma} \\
        & p(x) - \gamma = \sigma(x) \\
        & \sigma(x) \in \Sigma[x]_{\leq d}.
    \end{aligned}
\end{equation}

Equation \ref{eq:poly_sos} is a semidefinite program in terms of the Gram matrix $Q$ that defines $\sigma(x)$. The coefficient matching conditions of $p(x)-\gamma = \sigma(x)$ are linear constraints in the entries of $Q$. Given an optimal $Q$, the SOS decomposition of $\sigma(x)$ can be recovered by a matrix factorization of $Q$ (\emph{e.g.} Cholesky). The sequence $\{p_d^*\}$ is a increasing set of lower bounds to $p^*$.

% Sum-of-squares (SOS) methods offer a convex relaxation of polynomial nonnegativity: a polynomial $p(x) = \sum_{i} {p_i(x)^2} $ where $p_i(x) \in \R[x]$ is by definition nonnegative. Cones involved in this problem are the set of polynomials in $x$ with real coefficients $\R[x]$, polynomials of maximum total degree $d$ $\R[x]_{\leq d}$, and the subcone of sum-of-squares polynomials  $\Sigma$ $p(x) \in \Sigma_x$ has a Gram matrix description $p(x) = v(x)^\tr Q v(x)$ where $v(x)$ is a vector of monomials and $Q \succeq 0$. 

Constrained polynomial optimization can be approached through SOS methods. A basic semialgebraic set $\K$ is defined by a finite number of bounded-degree polynomial inequality and equality constraints:
\begin{equation}
    \K = \{ x \mid g_i(x) \geq 0, \quad h_j(x) = 0\}.
\end{equation} 
More generally, a semialgebraic set is the closure of basic semialgebraic sets under finite unions and projections down to coordinates. A constrained polynomial optimization problem can be formulated as Equation \ref{eq:poly_nonneg} with the added condition that $x \in \K$. If there exists a sufficiently large constant $R$ such that $\K \subset \{x \mid R - \norm{x}_2^2 \geq 0 \}$, then the
Putinar Positivstellensatz (Psatz) yields an equivalent formulation for polynomial constraint-multipliers $\zeta_i(x)$ and $\phi_j(x)$ \cite{putinar1993compact}:
% is a compact (Archimedean) semialgebraic set defined by inequalities $g_i$ and equalities $h_j$, optimization of $p(x)$ over $\K$ is equivalent to solving the following program ():
%
\begin{equation}
\label{eq:opt_K}
    \begin{aligned}
        p^* &= \max_{\gamma, \ \sigma, \ \zeta, \ \phi_j} {\gamma} \\
        & p(x) - \gamma = \sigma_(x) + \sum_i {\zeta_i(x)g_i(x)} + \sum_j{h_j(x) \phi_j(x)} \\
        & \sigma(x) \in \Sigma[x] \qquad \zeta(x) \in \Sigma[x] \qquad  \phi_j(x) \in \R[x].
    \end{aligned}
\end{equation}

%Bounding the degree of $p(x)$ (size of $Q$) forms an approximation hierarchy, which agrees with the global optimum if rank conditions are satisfied. Constrained optimization adds additional polynomial decision variables which may be relaxed to SOS and constrained in degree (Putinar multipliers).
%The Putinar expression $P(x) \geq0$ is made computationally tractable by relaxing $\sigma_i(x) \geq0 \rightarrow \sigma_i(x) \in SOS$ and by instituting a degree bound on polynomials $\sigma_i(x), \phi_j(x)$.

Equation \eqref{eq:opt_K} is an SDP when restricted to polynomials $\sigma(x), \zeta_i(x), \phi_j(x)$ of bounded degree $d$. 
If $\K$ satisfies an Archimedean condition then the sequence of lower bounds $p_d^* \leq p^*_{d+1} \leq \ldots$ will reach $p^*$ at a finite degree $d$\new{ (\new{Theorem 5.6 and 4.1 of} \cite{lasserre2010moments}).}
The size of the Gram matrix $Q$ scales as $O(N^{d})$ for $x \in \R^N$, and SDP performance is polynomial in $N$ \cite{lasserre2010moments}.

Utilizing sparsity in polynomial optimization can reduce computational complexity. Waki introduced a Correlative Sparsity Graph (CSP) $\gs(\vs, \es)$ where vertices $\vs$ are variables $x_i$ in the problem. An edge $(x_i, x_j) \in \es$ appears if $x_i$ and $x_j$ are multiplied together in a monomial in the cost $p(x)$, or if they appear together in at least one constraint $g_i(x)$ or $h_j(x)$ \cite{waki2006sums}. Let $\cs_k$ denote the cliques of the chordal-completed CSP graph, and $E_k x = x_k$ be the variables of $x$ present in CSP clique $\cs_k$. The Putinar Psatz in Equation \eqref{eq:opt_K} can formulated as a sparse problem over per-clique polynomials \cite{lasserre2006convergent}:
\begin{equation}
\label{eq:opt_K_sparse}
    \begin{aligned}
        p^* &= \max_{\gamma, \ \zeta_i, \ \phi_j} {\gamma} \\
        & p(x) - \gamma = \sum_k \sigma_{k}(x) + \sum_{i,k} {\zeta_{ik}(x)g_i(x)} + \sum_{j,k} {h_{jk}(x) \phi_j(x)}\\
        %&\sigma_i(x) \in SOS \qquad  \qquad  \phi_j(x) \in \R[x]
        &\zeta_{ik}(x_k) \in \Sigma[x_k] \qquad \sigma_k(x_k) \in \Sigma[x_k] \qquad \zeta_{ik}(x) \in \Sigma[x_k].
    \end{aligned}
\end{equation}

If the size of the largest CSP clique has cardinality $\kappa$, then the Psatz of degree $d$ in equation \eqref{eq:opt_K_sparse} will have expected $O(\kappa^{d})$ complexity. Other methods allow for decompositions according to the monomial structure in constraints. Josz uses monomial sparsity, which is a subset of the correlative sparsity graph \cite{josz2018lasserre}. Wang introduced a Term sparsity graph (TSSOS) based on links between monomials. The term sparsity graph allows for moment matrices to be block-diagonalized with an equivalent and mutually recoverable objective as SOS program \cite{wang2019tssos}. Term sparsity may be combined with Correlative Sparsity to maximize performance \cite{wang2021chordal}. Each of these methods use the structure of $p(x)$ and $\K$ to form computationally efficient semidefinite programs for polynomial optimization. 
\new{Another cone for sparse POP approximation as an alternative to SOS polynomials is the  Sum of Nonnegative Circuit polynomials \cite{iliman2016lower,  dressler2020global, wang2020second}. }

% The optimization problem in Equation \eqref{eq:opt_K} can be translated into a sparse Putinar form 

% \begin{equation}
% \label{eq:opt_K}
%     \begin{aligned}
%         p^* &= \max_{\gamma, \ \sigma_i, \ \phi_j} {\gamma} \\
%         & p(x) - \gamma = \sum_{i,k} {\sigma_{ik}(x)g_i(x)} + \sum_{j,k} {h_{jk}(x) \phi_j(x)} + \sum_k \zeta_{k}(x)\\
%         %&\sigma_i(x) \in SOS \qquad  \qquad  \phi_j(x) \in \R[x]
%         &\sigma_{ik}(x) \in \Sigma_{v_{ik}(x)} \\
%         &\zeta_k(x) \in \Sigma_{v_k(x)}
%         &\sigma_ik(x) \in \Sigma_{v_ik(x)}
%     \end{aligned}
% \end{equation}

% \citet{zheng2019bridge} leverage correlative sparsity to form optimization problems over the Sparse Sum of Squares cone of polynomials $SSOS(\es)$ ($Q \in \psd_+(\es,0)$).
%The general theory is that $p(x) = v(x)^\tr Q v(x)$ should  \citet{zheng2019bridge} chordally decompose the correlative sparsity graph to define a Sparse Sum of Squares cone of polynomials $SSOS(\es)$. 
%The correlative sparsity graph $\gs(\vs, \es)$ induces a sparse pattern $\gs_G(\vs_G, \es_G)$ on the Gram matrix $Q$, in which Agler's Theorem can be used for efficient optimization. More generally, the SDP realization of an SOS semialgebraic program induces an aggregate sparsity pattern on Gram matrices ${Q_i}_{i=0}^m \geq0$ involved in the problem ($P(x) = v(x)^\tr Q_0 v(x), \sigma_i(x) = v_i(x)^\tr Q_i v_i(x)$).

\subsection{Decomposed Structured Subsets for POPs}

Computational complexity can also be reduced by restricting polynomials to structured subsets. Given that $p(x) \in SOS$ implies that $p(x) = v(x)^\tr Q v(x), \ Q \geq 0$, setting $Q \in \dd$ or $Q \in \sdd$ results in polynomial cones $p(x) \in DSOS$ or $p(x) \in SDSOS$ respectively \cite{majumdar2014control}.

Decomposed structured subsets can be integrated into polynomial optimization. For a single SOS polynomial $p(x) \in SOS$, with chordal Gram sparsity pattern $\gs(\vs, \es)$ and maximal cliques $\{\cs_k\}_{k=1}^p$, Agler's theorem forms an equivalence of optima between $Q \in \psd_+^n$ and $Q \in \psd_+^n(\es,0)$. Restricting to standard structured subsets $p(x) \in DSOS$ forms $Q \in (\dd \cap \psd^n(\es,0)) = \dd^n(\es,0)$, and likewise $p(x) \in SDSOS$ forms $Q \in ( \sdd \cap \psd^n(\es,0)) = \sdd^n(\es,0)$. Allowing for mixed clique cones yields $\ks(\es,0)$, which may be broader than $\dd(\es,0)$ or $\sdd(\es,0)$ alone.

An example of decomposed structured subsets for polynomial optimization is the minimization of the (Rosenbrock-inspired) polynomial $f(x) = f_Q(x) + f_R(x)$, where 
\begin{align*}
   % f(x) &= f_R(x) + f_Q(x) \\
    f_R(x) &= \sum_{i=1}^{N-3} {10(x_{i+2} + 2 x_{i+1} - x_i^2)^2 + (1-x_i -x_{i+3})^2}, \\
    f_Q(x) &= x_{1:N/6}^\tr A x_{1:N/6}.
\end{align*}
%
%In a 120 variable problem, $f_Q$ is a quadric that involves the first $N/6=20$ variables. In this case, $A$ is the Lerner matrix defined as $A_{ij} = \text{min}(i/j, j/i)$. 
$f_Q(x)$ is a quadric where $A$ is the Lerner matrix defined as $\text{min}(i/j, j/i)$. 
With an 120 variable problem, the correlative sparsity graph of the Lehmer-Rosenbrock (LR) function $f(x)$ has cliques with size $\{4, 8, 11, 15, 231\}$. The size 15 occurs 97 times and all other clique sizes appear once. 
% The clique frequency is displayed in Figure \ref{fig:LR_CSP}.
The Lehmer-Rosenbrock function is chosen to highlight a set of SDPs with one very large clique. The standard Rosenbrock function is highly sparse, and does not have one giant component. 

% \begin{figure}[ht]
%     \centering
%     \includegraphics[width=0.5\linewidth]{img/LR_CSP_cliques_2.pdf}
%     \caption{\label{fig:LR_CSP}  Cliques for LR Correlative Sparsity Pattern}
% \end{figure}

% Performing this optimization with through a clique-decomposed TSSOS results in cliques of size $\{1\}$
Table \ref{tab:LR_CSP} shows the results of this optimization for the correlative sparsity pattern based on Sparse Sum of Squares \cite{zheng2019bridge}. This method is based on the sparsity of the Gram matrix $Q$ and optimizes over the cone of SOS polynomials of bounded degree. Table \ref{tab:LR_CSP} therefore shows a hierarchy of lower bounds to the minimum of $f(x)$ on $\R^{120}$. In the Cost and Time section, $K(\es, 0)$ is the cone where all cliques are in $K$, and $\ks(\es, 0)$ is the cone where the largest (231-sized) clique is in $K$ and all other cliques are in $\psd_+$. Upper bounds in this context would approximate the second-order moment relaxation, and would not give any useful bounds on the true polynomial optimum. Cells are merged if the cones are equal (as implemented).

\begin{table}[ht]
 \setlength{\abovecaptionskip}{0mm}
   \setlength{\belowcaptionskip}{0mm}
    \renewcommand\arraystretch{1.0}
\caption{\label{tab:LR_CSP} CSP $f(x)$ Lower Bounds }
\centering
\begin{tabular}{lcccc}
\toprule
         & \multicolumn{2}{l}{\hfil Cost}    & \multicolumn{2}{l}{\hfil Time (s)} \\
Cone     & $K^*(\es, 0)$            & $\ks^*(\es, 0)$      & $K^*(\es, 0)$            & $\ks^*(\es, 0)$        \\
\cmidrule(lr){2-3}\cmidrule(lr){4-5}
$\dd$    & -Inf.         & -Inf.         & 0.75         & 0.86           \\
$\sdd$   & -113.91         & -111.25        & 42.8        & 33.2      \\
$B_2$    & -111.63         & -111.19       & 15.7         & 12.6          \\
$B_3$    & -111.52       & -111.18       & 12.5         & 12.0          \\

$B_5$    & -111.06       & -111.05       & 12.75         & 12.6          \\
$B_{11}$   & \multicolumn{2}{l}{\hfil-110.85} & \multicolumn{2}{l}{\hfil27.1}    \\
$B_{21}$   & \multicolumn{2}{l}{\hfil-110.39} & \multicolumn{2}{l}{\hfil45.9}    \\
$B_{40}$   & \multicolumn{2}{l}{\hfil-110.26} & \multicolumn{2}{l}{\hfil112.7}    \\
$\psd_+$ & \multicolumn{2}{l}{\hfil-110.20} & \multicolumn{2}{l}{\hfil219.6}    \\
\bottomrule
\end{tabular}
\end{table}

\begin{figure}[ht]
    \centering
    \includegraphics[width=0.5\linewidth]{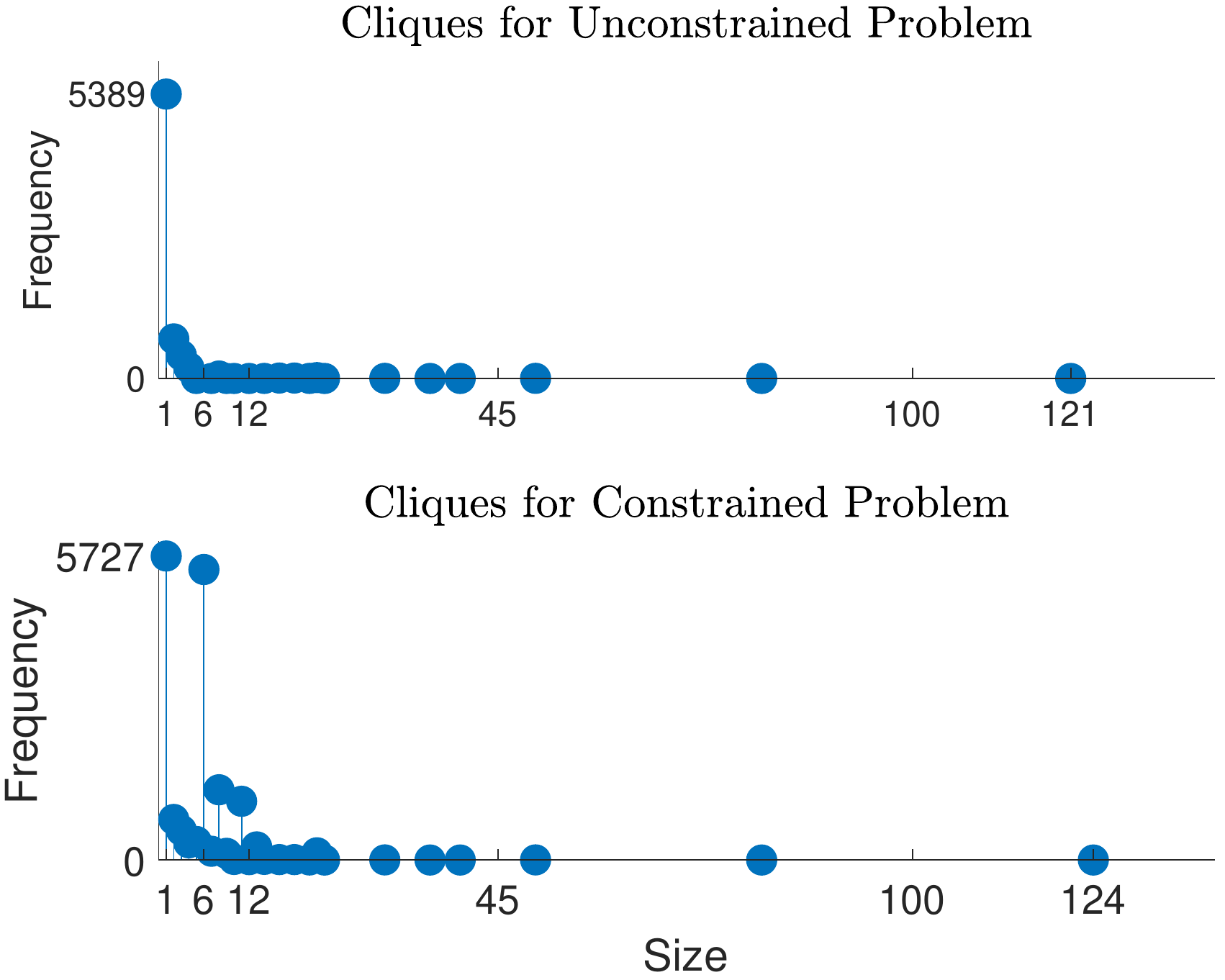}
    \caption{\label{fig:LR_cons} Cliques for LR TSSOS unconstrained and constrained optimization.}
\end{figure}

Constrained polynomial optimization offers additional freedom of cone-selection for decomposed structured subsets. Cones in $\ks$ can be chosen for each SOS constrained polynomial $p(x)$ and $\sigma_i(x)$. The following examples cover minimizing $f(x)$ over the semialgebraic set $x_i \in [1, 2]$. This region can be represented as $\K = \{g_i(x)=(x-1)(2-x) \geq 0 \}$. TSSOS \cite{wang2019tssos} supports constrained polynomial optimization problems and exploits term sparsity in the constraint functions $g(x)$. The TSSOS clique sizes of minimizing the LR function $f(x)$ over $\R^{120}$ (unconstrained) and over $\K$ (constrained) in the $d=2$ level of the Lasserre hierarchy are displayed in Figure \ref{fig:LR_cons}. This decomposition was obtained after the block hierarchy stabilized using the `clique' option in the TSSOS Julia implementation. TSSOS resulted in smaller clique sizes (compare 121 for TSSOS vs. 231 for CSP) as shown in Figure \ref{fig:LR_cons} at the expense of preprocessing time.

Figures \ref{fig:tssos_uncons} and \ref{fig:tssos_cons} display the time taken to provide lower bounds of $f(x)$ over the regions $\R^{120}$ and $[1,2]^{120}$ at the second Moment relaxation. Time is labeled in minutes to best provide contrast and intuition. Only values achieving the SDP optimum are displayed. In the constrained case, the lower bound of 4939.1 is attained in 31.8 minutes with the cone $B_3^*$ and size threshold 45, compared to 63.4 minutes on the full SDP $\psd_+(\es, ?)$.

\begin{table}[!ht]
    \centering
        \caption{Time to find unconstrained LR lower bound by TSSOS over $K^*$  (minutes)}
    \label{fig:tssos_uncons}
    \includegraphics[width=0.9\linewidth]{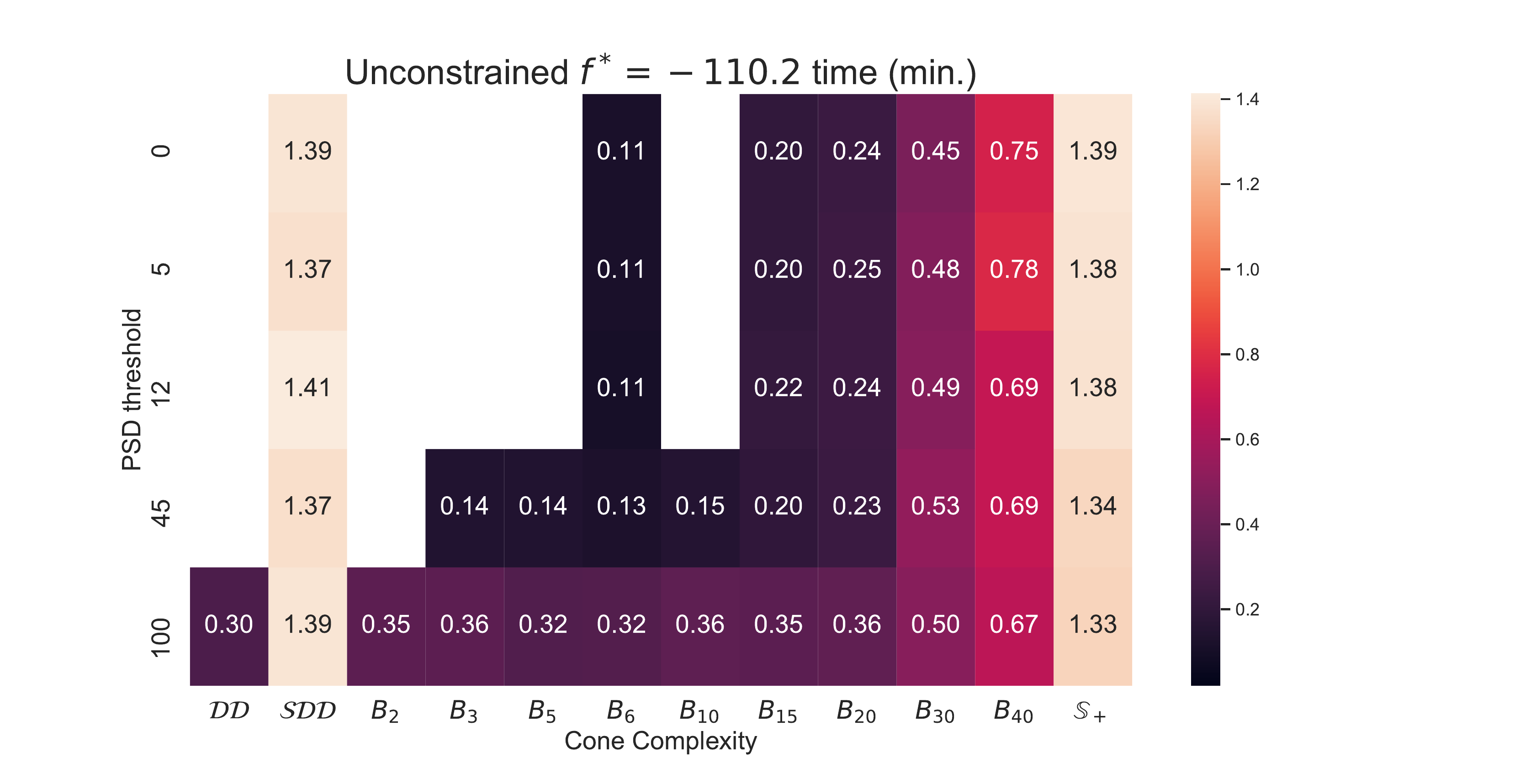}
\end{table}

\begin{table}[!ht]
    \centering
        \caption{Time to find constrained LR lower bound by TSSOS over $K^*$ (minutes)}
    \label{fig:tssos_cons}
    \includegraphics[width=0.9\linewidth]{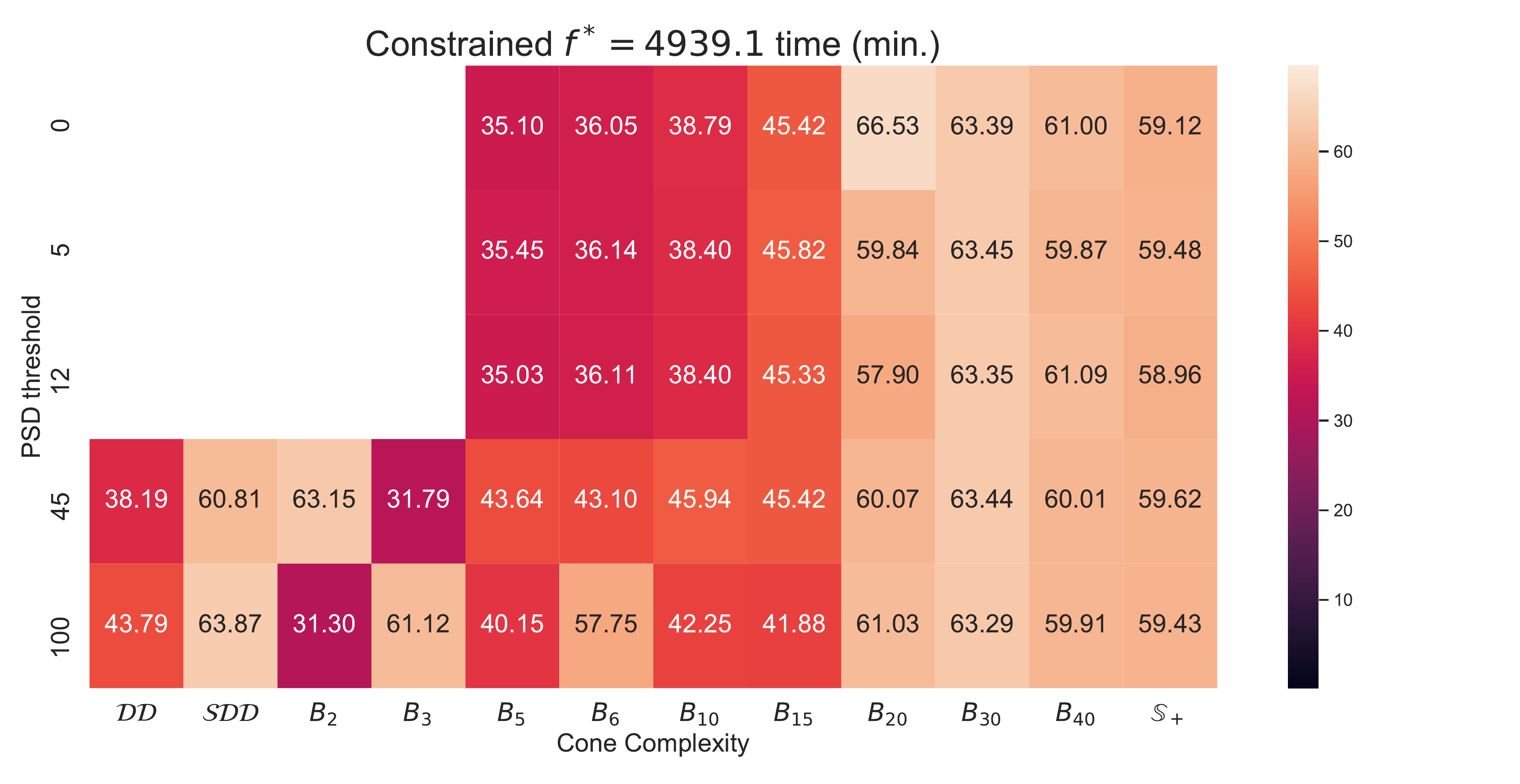}
\end{table}

\section{Conclusions} \label{sect:conclusion}
Structured subsets can be used to find upper and lower bounds of SDP optima. Decomposition methods may be able to convert large PSD constraints into smaller PSD blocks. This paper combines the two methods into decomposed structured subsets. Properties of these subsets are analyzed with their bound quality, and the facility to mix cones adds additional flexibility. Improved approximations are demonstrated on $\mathcal{H}_\infty$ norm and polynomial optimization problems. 
% Multiple forms of structure can be leveraged to find the tightest possible bounds.
Future directions include applying these techniques \new{to} network   $\mathcal{H}_{\infty}$-optimal control \new{and} more \new{POPs}. It would also be valuable to investigate compromises between cone complexity, additional consistency constraints, and approximation quality.
%Further extensions include solving truss and free material optimization problems with symmetry/*-algebra structure through decomposed structured subsets.
\newpage
\bibliographystyle{plain}
\bibliography{csos}

\newpage
\appendix

\section{Proof of Proposition 1}
\label{Sec:appendix}

In this Appendix, we provide the proof for Proposition 1 where $\es$ is not necessarily a chordal graph. %One direction is trivial

$\Leftarrow:$ If $Z_k \in \dd^{\vert \mathcal{C}_k \vert}$ then ${E_{\mathcal{C}_k}^\tr Z_k E_{\mathcal{C}_k}} \in \dd^n$ because diagonal elements of $Z_k$ will remain on the diagonal of ${E_{\mathcal{C}_k}^\tr Z_k E_{\mathcal{C}_k}}$, and the diagonal dominance relation is preserved in the embedded matrix.
% because the nonzero values of ${E_{\mathcal{C}_k}^\tr Z_k E_{\mathcal{C}_k}}$ are located only on the clique $\cs_k$. 
$\dd^n$ is a cone, so $ Z = \sum_{k=1}^p{E_{\mathcal{C}_k}^\tr Z_k E_{\mathcal{C}_k}} \in \dd^n$ for multiple cliques in $\cs$.
% $\Leftarrow:$ If $Z_k \in \dd^{\vert \mathcal{C}_k \vert}$, then ${E_{\mathcal{C}_k}^\tr Z_k E_{\mathcal{C}_k}} \in \dd^n$ and 
% $
%     % Z = \sum_{k=1}^p{E_{\mathcal{C}_k}^\tr Z_k E_{\mathcal{C}_k}} \in \dd^n.
% $
%
Let $\es^c \supseteq \es$ be a chordal completion of $\es$ with a clique cover $\cs^c$. All of the induced \new{submatrices of $Z$ by} cliques \new{in} $\cs^c$ are $DD$, so from Agler's theorem we have $Z \in \mathbb{S}^n(\mathcal{E},0)$.
% By Agler's theorem, we have $Z \in \mathbb{S}^n(\mathcal{E},0)$. 
Thus, $Z \in \dd^n(\mathcal{E},0)$.

$\Rightarrow:$ %Suppose $Z \in \dd^n(\mathcal{E},0)$. 
%All elements of a cone $K$ can be expressed as the convex combination of extreme points, which are points $x \in K$ such that there does not exist $x_1, x_2 \in K$ and $0 < \lambda  < 1$ such that $x = \lambda x_1 + (1-\lambda) x_2$. Extreme rays are the set $r x: r \geq 0$ for an extreme point $x$. Extreme rays of $\psd_+^n$ are the set of rank one symmetric matrices.
%
Let $e_i\in \mathbb{R}^{n}, e_j \in \mathbb{R}^{n}$ where $i \neq j$ be standard basis vectors, we define the following $DD$ basis matrices: 
%
%
%the set of symmetric diagonally dominant matrices can be generated by basis matrices 
% $$
% \begin{aligned}
% v^D_i &= e_i e_i^\tr , \\
% v^+_{ij}&=\begin{bmatrix}e_i & e_j\end{bmatrix}\begin{bmatrix}1 & 1 \\ 1 & 1 \end{bmatrix} \begin{bmatrix}e_i & e_j\end{bmatrix}^\tr , \\
% v^-_{ij}&=\begin{bmatrix}e_i & e_j\end{bmatrix}\begin{bmatrix}1 & -1 \\ -1 & 1 \end{bmatrix} \begin{bmatrix}e_i& e_j\end{bmatrix}^\tr.
% \end{aligned}
% $$
% $$
% \begin{aligned}
% v_i &= e_i e_i^\tr , \\
% v^+_{ij}&=(e_i + e_j)(e_i + e_j)^\tr , \\
% v^-_{ij}&=(e_i - e_j)(e_i - e_j)^\tr.
% \end{aligned}
% $$
%
\new{
\begin{align*}
   v_i &= e_i e_i^\tr, &
   v^\pm_{ij}&=(e_i \pm e_j)(e_i \pm e_j)^\tr. 
\end{align*}
}
Given a symmetric $Z \in \dd^n(\es, 0)$ for a (not-necessarily chordal) sparsity pattern $\es$, define the slack quantities $\Delta_i = Z_{ii} - \sum_{j\neq i}|Z_{ij}| \geq 0$. Such a $Z$ can be decomposed as:
\begin{align*}
Z &= \textstyle \sum_{i=1}^n \Delta_i v_i + \textstyle \sum_{(i,j) \in \mathcal{P} }Z_{ij}v^+_{ij} +  \textstyle \sum_{(i,j) \in \mathcal{N}}|Z_{ij}|v^-_{ij} \\
\mathcal{P} &:= \{(i,j) \mid  Z_{ij} >0, \  i < j\}\\
\mathcal{N} &:= \{(i,j) \mid  Z_{ij} <0, \ i < j\}. 
\end{align*}

% It is known that $\forall Z \in \dd^n$, there exist positive constants $c^D_i, c^+_{ij}, c^-_{ij}$ such that:
% \[
%     Z = \sum_i^{n} {c^D_i v^D_i} + \sum_{i \leq j}^n \left({c^+_{ij} v^+_{ij} + c^-_{ij} v^-_{ij}}\right).
% \]
% %
% %
% %An intuitive explanation of this is that every off-diagonal component is exactly offset by a positive diagonal component, and diagonal components can be larger than their required offsets.
% %
% Such a representation of $X$ with a minimal 
% $$
%     \sum_i^n{c_i^D} + \sum_{i \leq j}^n
%     \left(c^+_{ij} + c^-_{ij}\right)
% $$ 
% is unique. For any $Z_{ij} \neq 0$, then either $c^+_{ij} > 0$ or $c^-_{ij} > 0$. 
% Under this minimal representation, an entry where $Z_{ij} = 0$ will have $c^+_{ij} = 0$ and $c^-_{ij} = 0$. 

% When $Z \in \dd(\es, 0)$ for a (not-necessarily chordal) sparsity pattern $\es$, ${ij}$ corresponding to edges in $\es$ will likewise have either $c^+_{ij} > 0$ or $c^-_{ij} > 0$. 
By this characterization, $Z \in \dd^n(\es, 0)$ can be represented as the sum of $DD$ matrices with the same pattern $\mathcal{E}$ and clique cover $\cs$. %with coefficients $c_{ij}^+, c_{ij}^-$ corresponding to edges in $\es$. 
% By grouping the summands into different cliques, we get 
The terms $v_{ij}^+, \ v_{ij}^-$ with $Z_{ij} \neq 0$ can be uniquely assigned to some clique $(i, j) \in \cs_k$, and the slack terms $\Delta_i v_i$ can be distributed among all cliques in $\cs$ that include $i$. Grouping summands into cliques yields
\begin{align*}
Z &= \textstyle\sum_{k=1}^p{E_{\mathcal{C}_k}^\tr Z_k E_{\mathcal{C}_k}}, \quad 
Z_k \in \dd^{\vert \mathcal{C}_k \vert},
    & k=1,\,\ldots,\,p.\end{align*}
% The second statement is due to the fact that the sparsity pattern $\es$ of a matrix does not change after pre- or post-multiplication by a diagonal matrix.
Let $D$ be a positive definite (PD) diagonal matrix, and define matrices $D_{\cs_k} = E_{\cs_k} D E_{\cs_k}^\tr$ with inverses $D_{\cs_k}^{-1} = E_{\cs_k} D^{-1} E_{\cs_k}^\tr$. By definition \eqref{Eq:subsetdef} there exists a PD diagonal matrix $D$ for a $Z\in \sdd^n(\es, 0)$ such that $D Z D \in \dd^n$. Since pre and post-multiplying by a diagonal matrix does not change the sparsity pattern, $D Z D \in \dd^n(\es, 0)$. By the decomposition of $\dd^n$ matrices:
% A matrix $Z \in \sdd^n(\es, 0)$ must have an associated PD diagonal matrix $D$ such that $D Z D \in \dd^n(\es, 0)$ \eqref{Eq:subsetdef}, 
% From the definition of SDD matrices in \eqref{Eq:subsetdef}: $\forall Z \in \sdd^n(\es, 0)$, $\exists D \in \di^n_{>0}$ such that $D Z D \in \dd^n(\es, 0)$.
\begin{align*}
    D Z D = \textstyle \sum_{k=1}^p{E_{\mathcal{C}_k}^\tr \tilde{Z}_k E_{\mathcal{C}_k}},  \end{align*}
where $ \tilde{Z}_k \in \dd^{\abs{\cs_k}}, k = 1, \ldots, p$. This leads to
\begin{align*}
    % X &= D^{-1} \left(\sum_{k=1}^p{E_{\mathcal{C}_k}^\tr \tilde{X}_k E_{\mathcal{C}_k}}\right) D^{-1} \\
    Z = \textstyle \sum_{k=1}^p{D^{-1} E_{\mathcal{C}_k}^\tr \tilde{Z}_k E_{\mathcal{C}_k}D^{-1} } = \sum_{k=1}^p{ E_{\mathcal{C}_k}^\tr {Z}_k E_{\mathcal{C}_k} },
\end{align*}
where $Z_k = D^{-1}_{\cs_k} \tilde{Z}_k D^{-1}_{\cs_k} \in \sdd^{|\cs_k|}$, since $E_{\cs_k} D^{-1}= D^{-1}_{\cs_k} E_{\cs_k}$ and $D_{\cs_k}^{-1} = E_{\cs_k} D^{-1} E_{\cs_k}^\tr$, completing the proof.

\section{Combining Decompositions}
\label{sec:symm_appendix}
This section of the appendix reviews symmetry and *-algebra structure, shows how these structures can be incorporated into the Decomposed Structured Subset framework, and demonstrates how exploiting multiple forms of structure results in finer approximations.

\subsection{Symmetry/*-algebra Decomposition}

\label{sec:algebra}

An additional form of structure occurs when all constraint and cost matrices $(C, A_i)$ can be simultaneously block diagonalized with a unitary matrix $P$:
%\begin{align*} 
\[
\inp{C}{X} = \inp{P^T C P}{\tilde{X}} = \sum_{k} {\inp{\tilde{C}_k}{\tilde{X}_k}}\\
\]
%\inp{A_i}{X} = \inp{P^T A_i P}{\tilde{X}} = \sum_{k} {\inp{\tilde{A_i}_k}{\tilde{X}_k}}
%\end{align*}
Application of $P$ effectively breaks the large PSD variable $X = P \tilde{X} P^T$ into a product of smaller PSD variables $\tilde{X_k}$, and the SDP in $(\tilde{C}_k, \tilde{A_i}_k)$ will have an equivalent optimum as $(C, A_i)$. This block diagonalization can occur if all matrices $(C, A_i)$ lie in a common *-algebra $\A$. *-algebras are  subsets of matrices that are closed under addition, products, and transposition. If $(C, A_i)$ are all symmetric of size $n$, the blocks are free symmetric matrices of size $n_i$ with multiplicity $m_i$. This is expressed in the Wedderburn decomposition \citep{wedderburn1934lectures}:
\[ \A \cong \bigoplus_{i=1}^p {I_{m_i} \otimes S ^{n_i}} \qquad  \inp{C}{X} = \sum_i m_i \inp{\tilde{C}_i}{\tilde{X}_i}  \]

The Wedderburn decomposition of a *-algebra, given its basis, can be calculated numerically by randomized linear algebra \citep{murota2010numerical}. Group-invariant SDPs are a specific case of *-algebra structure. The Wedderburn decomposition of $G$-invariant matrices is related to the isotypic decomposition, and if $G$ is known in advance, the block-diagonalizing matrix $P$ may be calculated explicitly. Jordan decomposition generalizes symmetry and *-algebra decomposition, and involves projections onto an invariant subspace followed by splitting into the product of Jordan subalgebras (also requiring randomized computation) \citep{permenter2016dimension}. While preprocessing group and *-algebra structure is usually more expensive than exploiting chordal decompositions, their block-diagonalizability ensures that there are no overlaps between blocks (unlike clique-consistency constraints). 

If the semialgebraic set $\K$ is invariant under a symmetry group $G$, this structure can be exploited to reduce the computational burden. There exists an extensive literature on polynomial invariant and equivariant theory \citep{sturmfels2008algorithms}, and the main reference for SOS-based symmetry in polynomial optimization is \citep{gatermann2004symmetry}. The Hironaka decomposition of the invariant ring $\R[x]^G$ effects a block-diagonalization, which applies to $\K$ only if each polynomial $g_i(x), h_j(x)$ is invariant under $G$-actions. The methods of this paper can be applied after symmetry processing to restrict each block to a structured subsets, after first applying sparsity.

\subsection{Symmetry structure}

Decomposed structured subsets can be applied to symmetric SDPs in the same manner as chordally sparse problems. In the *-algebra framework, the cone-set $\ks$ refers to the cone of each symmetric block in the program.

\begin{figure}[ht]
    \centering
        \includegraphics[width=0.7\linewidth]{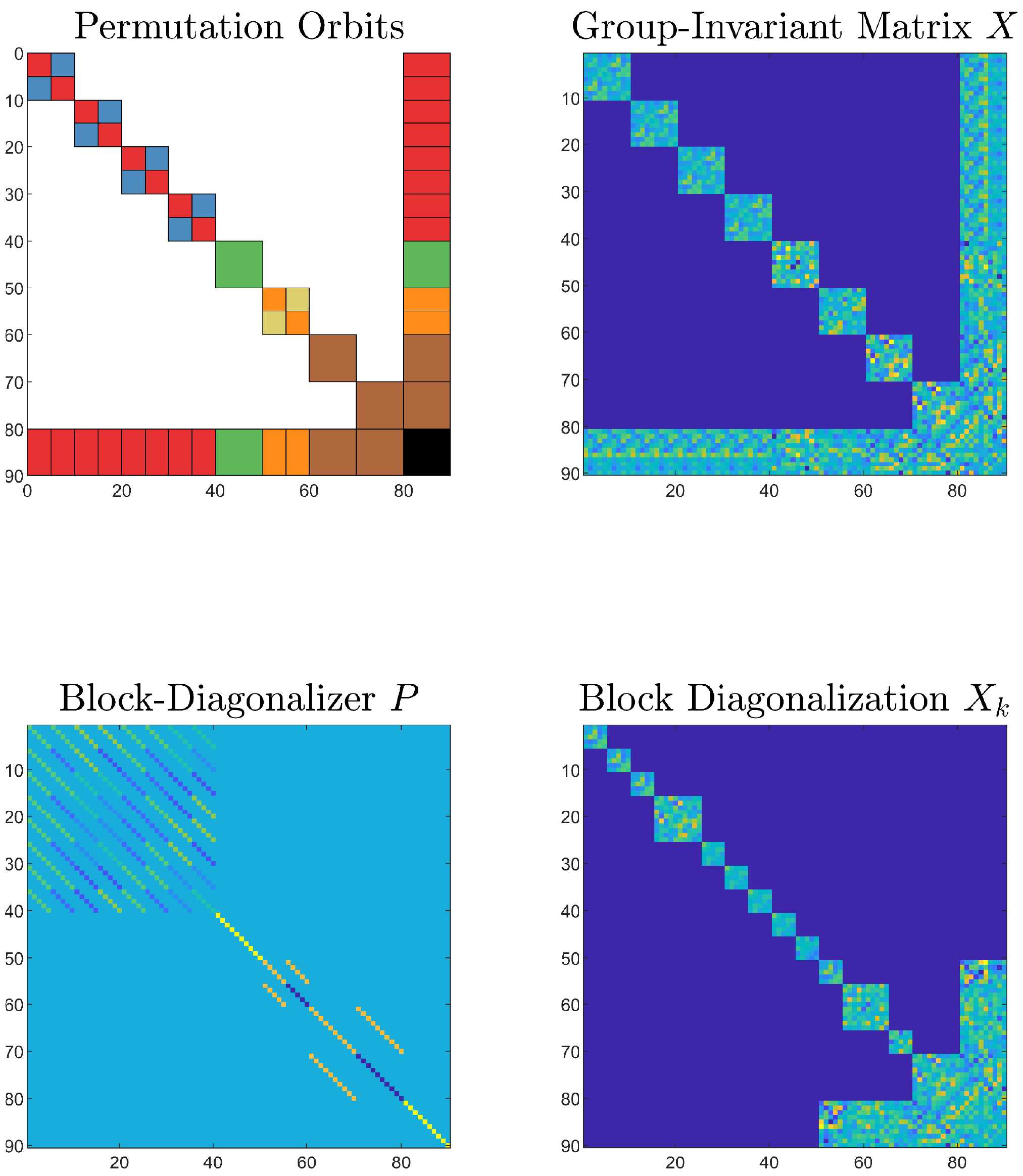}
    \caption{Sparse and Symmetric *-algebra}
    \label{fig:sym_split}
\end{figure}

A semidefinite program may contain more than one kind of structure. Figure \ref{fig:sym_split} shows block-arrow (sparse) matrices in the invariant ring $[\psd^{90}]^G$ under block-permutation action.
% where the group $G = S_4 \times 0 \times 0 \times \Z_2$ and $0$ is the one-element Zero group (identity action). 
The permutation structure is illustrated in the top-left pane of Figure \ref{fig:sym_split}. Matrices $X \in [\psd^{90}]^G$ are invariant under swapping blocks of the same color, and some blocks are additionally invariant under swapping the top and bottom halves (checkered pattern, split blocks). The permutation group $G$ acting on these matrices is:
\[G = (S_4 \times \Z_2) \times Id \times \Z_2 \times \Z_2\]

The top right pane of Figure \ref{fig:sym_split} shows a one such matrix $X \in [\psd^{90}]^G$. All such matrices can be block diagonalized (bottom right) under unitary action by the matrix $P$ (bottom left, based on the Discrete Cosine Transform). 
The block-sizes and multiplicities are:

\[(m,n) = \{(3, 5), (1, 10), (4,5), (1,5), (1, 40)^{sparse}\}\]
% Each orbit in the permutation has an component The block-orbits are of size $O = (4,1,1,2)$, so additional blocks appear outside the block-arrow section with multiplicities $O-1=(3, 0, 0, 1)$. 
% The Wedderburn decomposition has parameters:

% \[(m,n) = \{(3, 10), (1, 10), (1, 50)^{sparse}\}\]

% \begin{figure}[ht]
%     \centering
%     \includegraphics[width=0.5\linewidth]{img/sym_block_arrow.png}
%     \caption{Permutation invariant across blocks}
%     \label{fig:sym_block}
% \end{figure}

% Table \ref{tab:sym_sdp} shows the cost and time of running a randomly generated $G$-invariant SDP with 80 equality constraints under the cone $\sdd$. The columns show if the matrix has been block-diagonalized before applying $\sdd$, and the columns indicate if Grone's theorem has been applied to apply $\sdd$ to cliques.  (Full, Sparse) has the same cost as (Sym, Sparse), but takes longer to run because block-diagonalization eliminates redundancies in blocks with multiplicities.
%\cmidrule(lr){1-3}\cmidrule(lr){4-6}

%     \begin{table}[h]
%     \centering
% \resizebox{0.5\linewidth}{!} {\begin{tabular}{l|cccl|cc}
% Cost &   Full & Sym. & \qquad & Time (s) &   Full & Sym.\\ 
% \cmidrule(lr){1-3} \cmidrule(lr){5-7}
% Full & 10.82 & 10.37 & &Full & 55.07 & 12.39 \\
% Sparse &7.40 & 7.40 & & Sparse &13.39 & 7.11
% \end{tabular}}
%     \caption{\label{tab:sym_sdp} $\sdd$  Block Arrow with Symmetry (Fig. \ref{fig:sym_block})}
    
%     \end{table}

% Figure \ref{fig:sym_split} adds an additional symmetry in some of the blocks. 
Table \ref{tab:sym_sdp_split} shows the cost and time of running a randomly generated $G$-invariant SDP with 80 equality constraints under the cone $\sdd$. The columns show if the matrix has been block-diagonalized before applying $\sdd$, and the rows indicate if Grone's theorem has been applied to apply $\sdd$ to cliques.  (Full, Sparse) has the same cost as (Sym, Sparse), but takes longer to run because block-diagonalization eliminates redundancies in blocks with multiplicities.
Wedderburn Decomposition, whole blocks (green, brown orbits) contribute blocks of size 10 while split blocks have block-sizes of 5.

% Table \ref{tab:sym_sdp_split} displays similar results for a random $\tilde{G}$-invariant SDP. (Full, Sparse) has a higher cost and runtime than (Sym, Sparse). 

        \begin{table}[!ht]
            \caption{\label{tab:sym_sdp_split} $\sdd$  Block Arrow with Symmetry}
    \centering
 {
    \centering
    \begin{tabular}{l|cccl|cc}
Cost &   Full & Sym. & \qquad & Time (s) &   Full & Sym.\\ 
\cmidrule(lr){1-3} \cmidrule(lr){5-7}
Full & 12.96  &  10.86 & &Full & 124.5 & 19.3   \\
Sparse &9.49 & 8.44 & & Sparse & 38.2 & 12.0  
\end{tabular}}
    \end{table}

% \input{paper_sections_extended/symmetry_appendix}

% \cite{boyd1994linear}
% \maketitle
% \thispagestyle{empty}

\end{document}